\documentclass[12pt,a4paper]{article}%
\usepackage{amsmath}
\usepackage{epsfig}%
\usepackage{amsfonts}%
\usepackage{amssymb}
\usepackage{enumitem}
\usepackage{float}
\usepackage{lmodern}
\usepackage{mathtools}
\usepackage{tikz}
\usepackage{graphics}
\usepackage{graphicx}
\usepackage{cleveref}
\usepackage{multirow}
\usepackage{url}
\usepackage{pgfplots}
\usepackage{algorithm}
\usepackage[noend]{algpseudocode}
\def\BState{\State\hskip-\ALG@thistlm}

\makeatletter
\newcommand{\StatexIndent}[1][3]{%
  \setlength\@tempdima{\algorithmicindent}%
  \Statex\hskip\dimexpr#1\@tempdima\relax}
\usepackage{color}

    \usepackage[numbers]{natbib}
\usepackage[%
    font={small,sf},
    labelfont=bf,
    format=hang,    
    format=plain,
    margin=0pt,
    width=0.8\textwidth,
]{caption}
\tikzset{
    state/.style={
           rectangle,
           rounded corners,
           draw=black, very thick,
           minimum height=2em,
           inner sep=2pt,
           text centered,
           },
}
\usepackage{subcaption}

\usepackage{graphicx}
\usepackage[normalem]{ulem}
 \usepackage{authblk}
 
\newtheorem{theorem}{Theorem}[section]

\newtheorem{definition}[theorem]{Definition}

\crefname{lemma}{Lemma}{Lemmas}

\newtheorem{remark}{Remark}
\crefname{remark}{Remark}{Remarks}

\makeatletter\@addtoreset{equation}{section}\makeatother
\makeatletter\@addtoreset{figure}{section}\makeatother
\makeatletter\@addtoreset{table}{section}\makeatother
\newenvironment{customlegend}[1][]{%
\begingroup
\csname pgfplots@init@cleared@structures\endcsname
\pgfplotsset{#1}
}{%
\csname pgfplots@createlegend\endcsname
\endgroup
}%
\definecolor{background-color}{gray}{0.98}
\def\addlegendimage{\csname pgfplots@addlegendimage\endcsname}

\begin{document}

\title{Reduced basis method for the nonlinear Poisson-Boltzmann equation regularized by the 
range-separated canonical tensor format}
 
\author[1,2]{Cleophas Kweyu \thanks{\tt kweyuc@mu.ac.ke}}
\author[1]{Lihong Feng \thanks{\tt feng@mpi-magdeburg.mpg.de}}
\author[1]{Matthias Stein \thanks{\tt matthias.stein@mpi-magdeburg.mpg.de}}
\author[1,3]{Peter Benner \thanks{\tt benner@mpi-magdeburg.mpg.de}}
\affil[1]{Max Planck Institute for Dynamics of Complex Technical Systems, Sandtorstr.~1, D-39106 
Magdeburg, Germany}  
\affil[2]{Moi University, Department of Mathematics and Physics, P.O. Box 3900-30100, Eldoret, Kenya} 
\affil[3]{Otto von Guericke University, Faculty of Mathematics, Magdeburg, Germany}
  \date{}

\maketitle

\begin{abstract}
The Poisson-Boltzmann equation (PBE) is a fundamental implicit solvent continuum model for 
calculating the electrostatic potential of large ionic solvated biomolecules. However, its numerical 
solution encounters severe challenges arising from its strong singularity and nonlinearity. In 
\cite{BeKKKS:18, KwKKMB:18}, the effect of strong singularities was eliminated by applying the 
range-separated (RS) canonical tensor format \cite{BKK_RS:18, khor-DiracRS:2018} to construct a solution 
decomposition scheme for the PBE. The RS tensor format allows to derive a smooth approximation to the 
Dirac delta distribution in order to obtain a regularized PBE (RPBE) model. However, solving the RPBE is 
still computationally demanding due to its high dimension $\mathcal{N}$, where $\mathcal{N}$ is 
always in the millions. In this study, we propose to apply the reduced basis method (RBM) and 
the (discrete) empirical interpolation method ((D)EIM) to the RPBE in order to construct a reduced order 
model (ROM) of low dimension $N \ll \mathcal{N}$, whose solution accurately approximates the 
nonlinear RPBE. The long-range potential can be obtained by lifting the ROM solution back to the 
$\mathcal{N}$-space while the short-range potential is directly precomputed analytically, 
thanks to the RS tensor format. The sum of both provides the total electrostatic potential. The main 
computational benefit is the avoidance of computing the numerical approximation of the singular electrostatic 
potential. We demonstrate in the numerical experiments, the accuracy and efficacy of the reduced basis (RB) 
approximation to the nonlinear RPBE (NRPBE) solution and the corresponding computational savings over the 
classical nonlinear PBE (NPBE) as well as over the RBM being applied to the classical NPBE.
\end{abstract}



\section{Introduction}
\label{sec:Intro}

\begin{figure}[b]
  \centering
\captionsetup{width=\linewidth}
    \includegraphics[width=9.0cm]{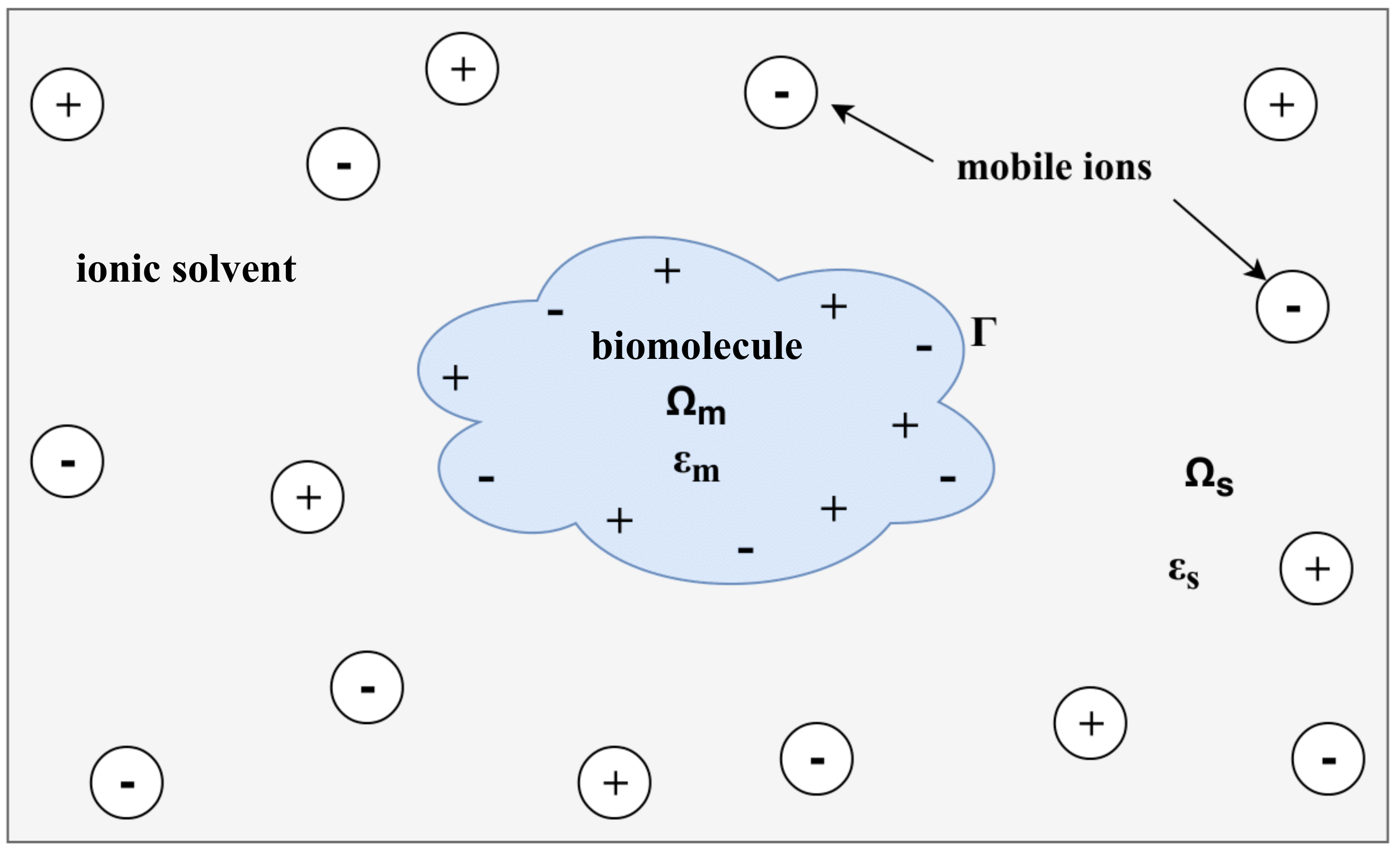}
    \caption{\label{fig:Biomolecular_system}
    2D representation of a biomolecular system.}
\end{figure}

The Poisson-Boltzmann equation (PBE) is a second order nonlinear elliptic partial differential 
equation (PDE) which is ubiquitous in the modeling of biochemical processes 
\cite{Baker:04,Baker2005,WaBa:04,Holst:2008}. It is used to calculate 
the electrostatic potential throughout the biomolecular system consisting of the biomolecule and the 
surrounding ionic or salt solution. More information about the significance of the electrostatic 
interactions and the related PBE post-processing, for instance, the electrostatic forces and energies, 
can be found in \cite{Fogolari2002, Neves-Petersen2003, Stein:2010}. \Cref{fig:Biomolecular_system} 
illustrates the two dimensional (2D) view of the biomolecular system consisting of a low dielectric 
molecular region ($\Omega_m$) encapsulated by an ionic solution of high dielectric in $\Omega_s$.

The analytical solution to the PBE for biomolecules with complex geometries, strong nonlinearities, and 
highly singular charge density distributions is not available \cite{Holst94, Dong2008}. To this end, 
numerical methods, for example, the finite difference method (FDM) \cite{Baker2001, Wang2010}, the finite 
element method (FEM) \cite{Baker2001, Holst2000}, the boundary element method (BEM) 
\cite{Boschitsch2004, Zhou1993}, are widely used to solve the PBE. Interested readers are referred to 
\cite{Holst:2008} for a thorough review of the aforementioned techniques for solving the PBE. However, 
the numerical solution to the PBE faces various challenges. The most severe are: the strong 
singularities, caused by the Dirac delta distribution sources; the strong nonlinearity, caused by the 
exponential nonlinear terms; the unbounded domain, due to the slow polynomial decay in the form of 
$1/\|\bar{x}\|$ as $\|\bar{x}\| \to \infty$; and imposing the correct jump or interface conditions to 
the rather irregular molecular domain, $\Gamma$ \cite{Xie:14, Chen:07, Mirzadeh:13}. 

In \cite{BeKKKS:18,KwKKMB:18}, the strong singularities are circumvented by using the range-separated (RS) 
canonical tensor format, which was introduced and analyzed in \cite{BKK_RS:18}. Consequently, 
a nonlinear regularized PBE (NRPBE) model is determined, which only solves for smooth long-range electrostatic 
potential. The jump conditions are annihilated due to the accurate splitting of the long- and short-range 
components of the total electrostatic potential using the RS tensor format. Nevertheless, the computational 
cost of solving the regularized nonlinear PBE is still high due to its high dimension 
$\mathcal{N} \approx \mathcal{O}(10^6)$. In this work, we apply the reduced basis method (RBM), in order to 
construct an accurate reduced order model (ROM) of much lower dimension, i.e., of $\mathcal{O}(10)$ for the 
NRPBE. The simulations for varying parameter values, in this case, the varying ionic strength, can be 
computed much faster by using the parametrized ROM and hence the Brownian dynamics simulations can be 
significantly accelerated. 

The RBM can also be directly applied to the classical NPBE without any regularization. We demonstrate that 
applying the RBM to the NRPBE gives rise to more robust and efficient solution to the problem as compared to
applying the RBM to the classical NPBE.

It is worth noting that the RBM has been applied to a simplified variant of the classical 
nonlinear PBE in \cite{Ji:2018} in $1$ and $2$ dimensions wherein smooth exponential functions were used as 
the source terms. In this work, we apply the RBM to the $3$-dimensional NRPBE for biomolecular simulation of 
large (complex) solvated biomolecules, for example, a protein in an ionic solutions, whose electrostatic 
potential $u(\bar{x})$ is characterized by the slow polynomial decay in $1/\|\bar{x}\|$, i.e., $u(\infty) = 0$, 
hence large domains have to be considered for accurate approximation of boundary conditions, resulting in 
high dimension of $\mathcal{O}(10^6)$ for the discretized system \cite{KwKKMB:18}. Furthermore, we consider 
Dirichlet boundary conditions which are nonaffine in parameter, thereby requiring the application of 
(D)EIM in order to reduce the complexity of the boundary conditions, 
and in turn, to further reduce the ROM complexity \cite{morChaS10,morGreMNetal07,morBarMNetal04}. 

The main contributions of this paper include: we have applied a more efficient numerical method for solving 
the NRPBE, which is based on first linearization via the Taylor series truncation of the nonlinear term, 
followed by discretization. This approach avoids the computation of the Jacobian of a huge matrix and also 
converges much faster than the standard Newton iteration; we have successfully applied the RBM and (D)EIM 
in order to reduce the dimension and complexity of the NRPBE; numerical comparison of RBM applied to the 
NRPBE and the classical NPBE shows that the former is more efficient and accurate.

The remainder of the paper is structured as follows. \Cref{sec:PBE_regularization} briefly reviews the 
approach of regularizing the PBE model by the RS canonical tensor format as proposed in 
\cite{BeKKKS:18,KwKKMB:18}. In \Cref{sec:RPBE_RBM}, the RBM framework and its application to both the 
regularized PBE (RPBE) and the classical PBE is introduced. Finally, \Cref{sec:Numer_Tests} presents the 
numerical experiments to illustrate the computational advantages of the RBM for the RPBE over the classical 
PBE. Comparisons with the solutions obtained by the standard FDM-based PBE solvers for the classical PBE 
are also presented.

\section{Mathematical model of the PBE}
\label{sec:PBE_model}

The nonlinear PBE for a symmetric 1:1 salt is given by 
\begin{equation}\label{eqn:PBE}
-\nabla\cdot(\epsilon(\bar{x})\nabla u(\bar{x})) + \bar{\kappa}^2(\bar{x})\sinh(u(\bar{x})) = 
  \sum_{i=1}^{N_m}q_i\delta(\bar{x}-\bar{x}_i),\quad \Omega \in \mathbb{R}^3,
\end{equation}
subject to 
\begin{equation}\label{eq:DH_solution}
u(\bar{x}) = \frac{1}{4\pi} \sum_{i=1}^{N_m}\frac{q_ie^{-\kappa(d-a_i)}}{\epsilon_s (1+\kappa a_i)d} 
\quad \mbox{on the boundary} \, \, \partial{\Omega}, \quad d = \lVert \bar{x}-\bar{x}_i \rVert, 
\quad \bar{x} = (x,y,z),      
\end{equation}
where $q_i = 4\pi e_c^2z_i/\kappa_B T$, $z_i$ is the partial charge of each atom, $u(\bar{x})$ represents 
the dimensionless potential, $\kappa_B T$, $\kappa_B$, $T$, and $e_c$ are the thermal energy, the Boltzmann 
constant, the absolute temperature, and the electron charge, respectively. The Debye-H\"uckel 
screening parameter, $\kappa^2 = {8\pi N_A e_c^2 I}/{1000\epsilon_s \kappa_BT}$, describes ion concentration 
and accessibility, $\epsilon_s$ is the solvent dielectric coefficient, $a_i$ is the atomic radius, and $N_m$ 
is the sum of the partial charges in the biomolecule. The sum of Dirac delta distributions represent the 
highly singular molecular charge density. 

The dielectric coefficient $\epsilon(\bar{x})$ and kappa function are piecewise constant functions 
given by
\begin{eqnarray}\label{eqn:diel_kappa_def}
 \epsilon(\bar{x}) =
  \begin{cases}
   \epsilon_m = 2 & \text{if } \bar{x} \in \Omega_m\\
   \epsilon_s \,\,= 78.54 & \text{if } \bar{x} \in \Omega_s
  \end{cases}, \quad \quad
 \bar{\kappa}(\bar{x}) =
  \begin{cases}
   0 & \text{if } \bar{x} \in \Omega_m \\
    \sqrt{\epsilon_s}\kappa & \text{if } \bar{x} \in \Omega_s 
  \end{cases},
\end{eqnarray}
where $\Omega_m$ and $\Omega_s$ are the regions occupied by the protein molecule and by the ionic 
solution, respectively, as shown in \Cref{fig:Biomolecular_system}. See \cite{SharpHonig90,Holst94, 
KwBFM2017} for discussions regarding the PBE theory and the importance of (\ref{eqn:PBE}) in biomolecular 
modeling.

The PBE in (\ref{eqn:PBE}) can be linearized for small electrostatic potentials by retaining the 
first term of the Taylor series expansion of the nonlinear function $\sinh(u(\bar{x}))$ \cite{Fogolari99}. 
The LPBE is thus given by
\begin{equation}\label{eq:LPBE}
   -\nabla\cdot(\epsilon(\bar{x})\nabla u(\bar{x})) + \bar{\kappa}^2(\bar{x})u(\bar{x}) = 
   \sum_{i=1}^{N_m}q_i\delta(\bar{x}-\bar{x}_i).
\end{equation}
The LPBE is much easier to solve and very accurate for lowly charged biomolecules, for example, 
proteins. However, for highly charged biomolecules, such as nucleic acids, it is not as 
accurate as the nonlinear variant due to the magnitude of the electric field at the interface 
between the solute and the solvent, $\Gamma$ \cite{Qin_Ray_2010, Fogolari2002}.

\section{Regularization of the PBE by the RS tensor format}
\label{sec:PBE_regularization}

The numerical approximation of the PBE is hindered by the highly singular sources described by a sum 
of Dirac delta distributions. This is because, for every singular charge $z_i$ in (\ref{eqn:PBE}), 
there corresponds degenerate behaviour in the electrostatic potential $u(\bar{x}_i)$ at each atomic 
position $\bar{x}_i$ in $\Omega_m$. To circumvent this drawback, various researchers have developed 
solution decomposition approaches for the PBE \cite{Xie:14,Mirzadeh:13,Chen:07,Chern:2003}. A common 
feature of these approaches is that they circumvent the building of numerical approximations corresponding 
to the Dirac delta distributions by solving a regularized PBE model for the smooth long-range electrostatic 
potential. This is enhanced by the fact that analytical expansions by the Newton kernel are possible in 
the solute sub-region $\Omega_m$.

In principle, the solution decomposition techniques for the PBE involve coupling of two equations 
for the electrostatic potential in the solute and solvent regions, through the interface, $\Gamma$ 
\cite{Chen:07, Chern:2003}. Due to the absence of ions within the molecular region $\Omega_m$, it is 
modeled by the Poisson equation,   
\begin{equation}\label{eqn:PE}
 -\nabla \cdot(\epsilon_m\nabla u) = \sum_{i=1}^{N_m}q_i \delta(\bar{x}-\bar{x}_i)  
 \quad \mbox{in} \,\, \Omega_m.
\end{equation}
On the other hand, no atoms are present in the solvent 
region $\Omega_s$, hence the charge density is purely modeled by the Boltzmann distribution, leading to
\begin{equation}\label{eqn:Boltzmann_distn}
 -\nabla \cdot(\epsilon_s\nabla u) + \bar{\kappa}^2\sinh(u) = 0  \quad \mbox{in} \,\, \Omega_s.
\end{equation}
Therefore, the two equations (\ref{eqn:PE}) and (\ref{eqn:Boltzmann_distn}) are coupled together 
via the jump (interface) boundary conditions
\begin{equation}\label{eqn:interface_condn}
 \left[ u\right]_{\Gamma} = 0, \quad \mbox{and} \quad \left[ \epsilon \frac{\partial u}{\partial 
 n_{\Gamma}}\right]_{\Gamma} = 0,
\end{equation}
where $\Gamma := \partial\Omega_m = \partial\Omega_s \cap \Omega_m$ and $\left[ f\right]_{\Gamma} = 
\lim\limits_{t \to 0} \left(f(\bar{x}+tn_{\Gamma}) - f(\bar{x}-tn_{\Gamma})\right)$. Here, we denote $n_{\Gamma}$ as
the unit outward normal direction of the interface $\Gamma$. 

In \cite{BeKKKS:18,KwKKMB:18}, the authors employ the RS canonical tensor format, developed and analyzed in 
\cite{BKK_RS:18}, to construct the solution decomposition of the PBE. This is realized 
by approximating the singular sources with a smooth function derived from the long-range component of 
the Newton potential sum. The resultant regularized PBE solves for the long-range electrostatic potential, 
which is then added to the short-range component that is precomputed from the RS tensor splitting of the 
Newton kernel. The regularized PBE (RPBE) model has demonstrated to be much more accurate than the 
classical PBE model in \cite{BeKKKS:18,KwKKMB:18}. We highlight the core ingredients for obtaining the 
RPBE in \Cref{ssec:Coulomb}.

\subsection{Canonical tensor representation of the Newton kernel}
 \label{ssec:Coulomb}

\begin{definition}
The \emph{Newton potential} of an integrable function (or a Radon measure) $f$ with compact support in 
$\mathbb{R}^3$ is defined as the convolution 
\begin{equation}\label{Newt_pot_compact}
u(\bar{x}) = \Gamma_N * f(\bar{x}) = \int_{\mathbb{R}^3}\Gamma_N(\bar{x} -\bar{y}) f(\bar{y})dy,
\end{equation}
where the Newton kernel $\Gamma_N = 1/\|\bar{x}\|$, has a mathematical singularity at the origin, and 
$\bar{y} \in \mathbb{R}^3$ \cite{Newt_kernel}. The Newton potential $u(\bar{x})$ satisfies the Poisson 
equation
\begin{equation}
 -\Delta u = f,
\end{equation}
where $f$ in this case is the source term of the system as defined in (\ref{eqn:Poisson_rhs_fn}). 
\end{definition}

Consider the single particle Newton potential (or the Newton kernel) $1/\|\bar{x}\|$, $\bar{x} \in 
\mathbb{R}^3$, which is a fundamental solution to the Poisson equation. It is well known that determining 
a weighted sum of interaction potentials (or Newton kernels), $P_N(\bar{x})$ in a large $N_m$-particle 
system, with the particle locations at $\bar{x}_i \in \mathbb{R}^3$, $i=1,...,N_m$, i.e.,
\begin{equation}\label{eqn:PotSum1}
 P_N(\bar{x})= \sum_{i=1}^{N_m} \frac{q_i}{\epsilon_m\|\bar{x}-\bar{x}_i\|}, \quad \bar{x}_i, 
 \bar{x} \in \Omega=[-b,b]^3,
\end{equation}
is quite computationally demanding. The Newton kernel exhibits a slow polynomial decay in $1/\|\bar{x}\|$ 
as $\|\bar{x}\| \to \infty$. Obviously, it has a singularity at $\bar{x}=(0,0,0)$, making its accurate 
grid representation difficult. The RS tensor format \cite{BKK_RS:18} can be exploited to construct an 
efficient grid-based technique for the calculation of $P_N(\bar{x})$ in multiparticle 
systems.

\begin{remark}\label{remk:Newton_PE}
 Notice that the Newton potential $P_N(\bar{x})$ in (\ref{eqn:PotSum1}) is a special case of 
 (\ref{Newt_pot_compact}) for a non-compact function 
 \begin{equation}\label{eqn:Poisson_rhs_fn}
 f(x) = \frac{1}{\epsilon_m}\sum_{i=1}^{N_m}q_i \delta(\bar{x}-\bar{x}_i).
 \end{equation}
\end{remark}

To obtain the canonical tensor representation of the Newton kernel, we follow the procedure in 
\cite{BKK_RS:18}, whereby we first consider the computational domain $\Omega=[-b,b]^3$, and introduce the 
uniform ($n^{\otimes 3}$) \footnote{$n^{\otimes 3} = n \times n \times n$ is a tensor representation of the 
3D Cartesian grid.} rectangular Cartesian grid $\Omega_{n}$ with mesh size $h=2b/n$ ($n$ even). 
Let $\{ \psi_i\}$ be a set of tensor-product piecewise constant basis functions, 
$\psi_i(\bar{x})=\prod_{\ell=1}^3 \psi_{i_\ell}^{(\ell)}(\bar{x}_\ell)$,
for the $3$-tuple index $i=(i_1,i_2,i_3)$, $i_\ell \in I_\ell=\{1,...,n\}$, $\ell=1,\, 2,\, 3 $.
The goal is to discretize the Newton kernel by its projection onto $\{ \psi_i\}$ as follows
\begin{equation}  \label{galten}
\mathbf{P}:=[p_i] \equiv [p(i_1,i_2,i_3)] \in \mathbb{R}^{n^{\otimes 3}}, \quad
 p_i = \int_{\mathbb{R}^3} \frac{\psi_i ({\bar{x}})}{\|\bar{x}\|} \,\, \mathrm{d}{\bar{x}},
\end{equation}
where $p_i$ is obtained from the vectors of the canonical tensor representation of the Newton kernel.

Next, determine the Laplace-Gauss transform representation of $1/\|\bar{x}\|$, and then apply the 
exponentially convergent sinc-quadrature approximation to obtain the separable expansion
\begin{equation}\label{eqn:sep_expansion}
    \frac{1}{\|\bar{x}\|} = \frac{2}{\sqrt{\pi}}\int_{\mathbb{R}^+} e^{-t^2\|\bar{x}\|^2}\mathrm{d}{t} 
    \approx  \sum\limits_{k=-M}^{M}a_k e^{-t_k^2\|\bar{x}\|^2} 
    = \sum\limits_{k=-M}^{M} a_k\prod\limits_{\ell=1}^{3} e^{-t_k^2\bar{x}_{\ell}^2},
  \end{equation}
where the quadrature points and weights in (\ref{eqn:sep_expansion}) are given by 
  \begin{equation}\label{eqn:quadrature_points}
    t_k = k \mathfrak{h}_M, \quad a_k = 2\mathfrak{h}_M/{\sqrt{\pi}}, \quad \mbox{with} \:
    \mathfrak{h}_M = C_0\log(M)/M, \quad C_0 \approx 3. 
  \end{equation}
  
\begin{figure}[t]
    \centering
    \captionsetup{width=\linewidth}
    \begin{subfigure}[b]{0.32\textwidth}
        \centering
        \includegraphics[height=3.5cm]{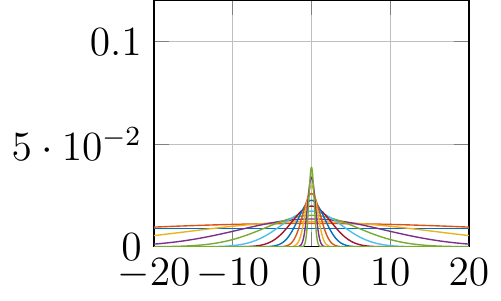}
        \caption{Long-range vectors.}
        \label{fig:Long-range_vecs}
    \end{subfigure}
    \hfill
    \begin{subfigure}[b]{0.32\textwidth}
        \centering
        \includegraphics[height=3.5cm]{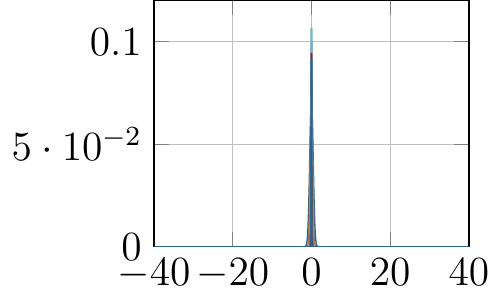}
        \caption{Short-range vectors.}
        \label{fig:Short-range_vecs}
    \end{subfigure}
    \hfill
    \begin{subfigure}[b]{0.32\textwidth}
        \centering
        \includegraphics[height=3.8cm]{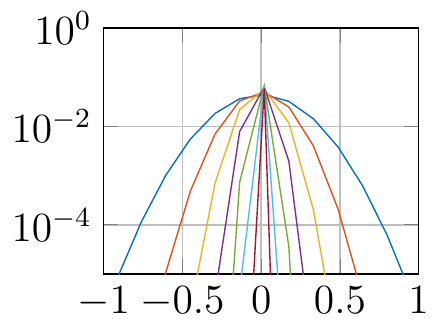}
        \caption{Log plot for \Cref{fig:Short-range_vecs}.}
        \label{fig:short_vecs_log}
    \end{subfigure}
    \caption{Canonical vectors for $n=1024$, $R=20$, and $R_l=12$.}
    \label{fig:Long_short_vecs}
\end{figure}

The mode three tensor $\textbf{P}$, can be approximated by the $R$-term $(R = 2M+1)$ canonical tensor
representation 
\begin{equation} \label{eqn:canon_repr}
\mathbf{P} \approx \mathbf{P}_R = \sum\limits_{k=-M}^{M} {\bf p}^{(1)}_k \otimes {\bf p}^{(2)}_k \otimes {\bf p}^{(3)}_k
\in \mathbb{R}^{n^{\otimes 3}},
\end{equation}
where ${\bf p}^{(\ell)}_k \in \mathbb{R}^n$ are obtained by substituting (\ref{eqn:sep_expansion}) 
into (\ref{galten}) and $``\otimes"$ \footnote{The outer product of two vectors, i.e., 
$x\otimes y = xy^T$ is a rank-one matrix, that of three vectors, i.e., $x\otimes y \otimes z$ is a 
rank-one tensor, and so forth.} denotes the outer (or tensor) product of vectors. For more details, see \cite{BKK_RS:18, KwKKMB:18}. Upon splitting the 
reference canonical tensor representation $\mathbf{P}_R$ by the procedure presented in \cite{BKK_RS:18}, 
we obtain the following decomposition
 \[
  \mathbf{P}_R = \mathbf{P}_{R_s} + \mathbf{P}_{R_l},
\]
where
\begin{equation} \label{eqn:Split_Tens}
    \mathbf{P}_{R_s} =
\sum\limits_{k\in {\cal K}_s} {\bf p}^{(1)}_k \otimes {\bf p}^{(2)}_k \otimes {\bf p}^{(3)}_k,
\quad \mathbf{P}_{R_l} =
\sum\limits_{k\in {\cal K}_l} {\bf p}^{(1)}_k \otimes {\bf p}^{(2)}_k \otimes {\bf p}^{(3)}_k.
\end{equation}
Here, ${\cal K}_l := \{k|k = 0,1, \ldots, R_l\}$ and ${\cal K}_s := \{k|k = R_l+1, \ldots, M\}$ are 
the sets of indices for the long- and short-range canonical vectors. The cross-sectional view of the 
respective localized and global vector components of the Newton potential in (\ref{eqn:Split_Tens}) on 
the $x$-axis is illustrated in \Cref{fig:Long_short_vecs}.

\newlength{\fwidth}
\newlength{\fheight}
\setlength{\fwidth}{3.9cm}
\setlength{\fheight}{3.7cm}

The aforementioned results are only valid for a single particle potential (or the Newton kernel, 
$1/\|\bar{x}\|$). In the case of a potential sum generated by a multiparticle system, for example, in 
(\ref{eqn:PotSum1}), the two components in (\ref{eqn:Split_Tens}) are treated independently due to their 
differences as far as their effective supports are concerned \cite{BKK_RS:18}. The following is an 
overview of the RS canonical tensor representation of $P_N(\bar{x})$ in (\ref{eqn:PotSum1}).

We first consider the tensor representation of only the long-range component ${\bf P}_{R_l} 
\in \mathbb{R}^{n^{\otimes 3}}$ which can be constructed by a direct sum of shift-and-windowing 
transforms, ${\cal W}_i$, of the reference tensor $\widetilde{\bf P}_{R_l} \in 
\mathbb{R}^{2n \times 2n \times 2n}$ from a large $(2n)^{\otimes 3}$ domain onto the original 
$n^{\otimes 3}$ domain. See \cite{VeBoKh:Ewald:14, BKK_RS:18} for detailed information.
\begin{equation}\label{eqn:Long-Range_Sum} 
{\bf P}_l = \sum_{i=1}^{N_m} {z_i} \, {\cal W}_i (\widetilde{\mathbf{P}}_{R_l})=
\sum_{i=1}^{N_m} {z_i} \, {\cal W}_i
(\sum\limits_{k\in {\cal K}_l} \widetilde{\bf p}^{(1)}_k \otimes \widetilde{\bf p}^{(2)}_k
\otimes \widetilde{\bf p}^{(3)}_k).
\end{equation}

\begin{remark}
 Note that ${\bf P}_l$ comprises of a collection of ${\bf P}_{R_l}$ at each atomic postion in 
 the entire protein, which have been shifted and windowed by the transform ${\cal W}_i$ of the 
 reference tensor $\widetilde{\bf P}_{R_l} \in \mathbb{R}^{2n \times 2n \times 2n}$. Clearly, 
 $\widetilde{\bf P}_{R_l}$ consists of ${\bf P}_{R_l}$ in a $(2n)^{\otimes 3}$ domain. 
\end{remark}

The reference tensor $\widetilde{\bf P}_{R_l}$ is mapped onto its sub-tensor of smaller size 
$n^{\otimes 3}$, by first shifting the center of $\widetilde{\bf P}_{R_l}$ to the grid-point $x_i$, 
and then windowing (restricting) the result onto the computational grid $\Omega_n$. The particle 
charges are denoted by $z_{i}$. The canonical rank of the tensor sum ${\bf P}_l$, of rank $RN_m$, was 
proven in \cite{BKK_RS:18} to depend only logarithmically on the number of particles $N_m$ involved in 
the summation. 

\begin{remark}
It is worth noting that for large biomolecules, the rank $RN_m$ and the $n^{\otimes3}$ Cartesian grid can 
be very large due to large $N_m$. In such cases, the canonical-to-Tucker (C2T) and the Tucker-to-canonical 
(T2C) transforms can be applied in order to obtain a low rank canonical tensor representation which 
accurately approximates the original tensor. The C2T transform employs the reduced higher order singular 
value decomposition (RHOSVD) to accomplish the rank reduction process \cite{khor-ml-2009}. 
\end{remark} 

On the other hand, the short-range part of the total electrostatic potential is represented by a 
single small size tensor ${\bf P}_s \in \mathbb{R}^{n^{\otimes 3}}$, known as the cumulated canonical 
tensors (CCT) \cite{BKK_RS:18}. The CCT comprises of localized subtensors whose effective supports are 
nonintersecting
\begin{equation}\label{eq:CCT_tensor}
 {\bf P}_s = \sum_{i=1}^{N_m} z_i {\bf U}_i, \quad \quad {\bf U}_i \in \mathbb{R}^{n_s^{\otimes 3}}, 
            \quad n_s \ll n,
\end{equation}
where $\mbox{diam}(\mbox{supp}{\bf U}_i) \leq 2\sigma_i$. Here, $\sigma_i$ is the atomic radius of 
each atom in the biomolecule.
\begin{remark}
 Notice that for biomolecules whose atoms have varying radii, we adjust the computation of the 
 short- and long-range range electrostatic potential accordingly by assigning the corresponding 
 vectors from \Cref{fig:Long_short_vecs} to atomic clusters of similar radii \cite{KwKKMB:18}.  
\end{remark}

\subsection{Construction of the nonlinear RPBE (NRPBE)} \label{sec:RPBE_RS}
We now have sufficient information to facilitate the construction of the NRPBE based on the simple 
splitting of the Dirac delta distribution \cite{khor-DiracRS:2018}. To fix the idea, from 
\Cref{remk:Newton_PE}, the weighted sum of interaction potentials in a large $N_m$-particle system as in 
(\ref{eqn:PotSum1}) is also the analytical solution to the Poisson equation (PE), i.e.,
\begin{equation}\label{eqn:PE_substit}
 -\epsilon_m\Delta P_N(\bar{x}) = \sum_{i=1}^{N_m}q_i \delta(\bar{x}-\bar{x}_i)  
 \quad \mbox{in} \,\, \mathbb{R}^3.
\end{equation}  

Consider the RS tensor splitting of the multiparticle Newton potential into a sum of long-range tensors 
${\bf P}_l$ in (\ref{eqn:Long-Range_Sum}) and a CCT tensor ${\bf P}_s$ in (\ref{eq:CCT_tensor}), i.e., 
\begin{equation}\label{eqn:Newt_splitting}
 \mathbf{P}_N(\bar{x}) = {\bf P}_s(\bar{x}) + {\bf P}_l(\bar{x}).
\end{equation}
Applying the discretized Laplacian operator to each component of $P_N(\bar{x})$, we obtain, 
\begin{equation}\label{eqn:Dirac_splitting}
 f^s:= -A_{\Delta} {\bf P}_s, \quad 
 \mbox{and} \quad f^l:= -A_{\Delta} {\bf P}_l,
\end{equation}
where $A_{\Delta}$ is the 3D finite difference Laplacian matrix defined on the uniform rectangular grid as
\begin{equation}\label{eqn:Lapl_Kron3}
A_{\Delta} = \Delta_{1} \otimes I_2\otimes I_3 + I_1 \otimes \Delta_{2} \otimes I_3 + 
I_1 \otimes I_2\otimes \Delta_{3},
\end{equation}
where $-\Delta_\ell = h_\ell^{-2} \mathrm{tridiag} \{ 1,-2,1 \} \in \mathbb{R}^{n_\ell \times n_\ell}$, 
$\ell=1,2,3$, denotes the discrete univariate Laplacian and $I_\ell$, $\ell=1,2,3$, is the identity matrix 
in each dimension. See \cite{BeKKKS:18, KwKKMB:18,khor-DiracRS:2018} for more details.

The nonlinear regularized PBE (NRPBE) can now be derived as follows. First, the unknown solution (or target 
electrostatic potential) $u$ to the PBE (\ref{eqn:PBE}) can be decomposed as $u = u^s + u^r$, where $u^s$ 
is the known singular function (or short-range component) and $u^l$ is the unknown long-range component to 
be determined. Therefore, the PBE (\ref{eqn:PBE}) can be rewritten as
\begin{equation}\label{eqn:PBE_splitting}
  \begin{rcases}
  \begin{aligned}
 -\nabla\cdot(\epsilon\nabla (u^s + u^r)) + \bar{\kappa}^2\sinh(u^s + u^r) &= f^s + f^l  
 \quad \mbox{in} \,\, \mathbb{R}^3,\\
 u &= g, \quad \mbox{on} \quad \partial\Omega,
  \end{aligned}
  \end{rcases}
\end{equation}
where the right-hand side of (\ref{eqn:PBE}) is replaced by $f^s + f^l $ due to (\ref{eqn:PE_substit}) and 
(\ref{eqn:Dirac_splitting}) and $g$ is the Dirichlet boundary conditions defined in (\ref{eq:DH_solution}).

It was proved and demonstrated in \cite{BeKKKS:18} that the function $f^s$ and the corresponding 
short-range potential $u^s$ are localized within the molecular region $\Omega_m$ and vanishes on the 
interface $\Gamma$. Moreover, in the PBE (\ref{eqn:PBE}), the function $\bar{\kappa}$ is piecewise constant 
as defined in (\ref{eqn:diel_kappa_def}), and $\bar{\kappa} = 0$ in $\Omega_m$. Therefore, we can rewrite 
the Boltzmann distribution term in (\ref{eqn:PBE_splitting}) as
 \begin{equation}\label{eqn:Boltzmann_lr}
  \bar{\kappa}^2\sinh(u^s + u^r) = \bar{\kappa}^2\sinh(u^r), \quad \mbox{because } \, u^s = 0 \,
  \mbox{ in } \Omega_s.
 \end{equation}

Consequently, following the splitting of the Dirac-delta distributions in (\ref{eqn:PBE_splitting}), the 
short-range component of the potential satisfies the Poisson equation, i.e.,
\begin{equation}\label{eqn:PE_sr}
 -\nabla\cdot(\epsilon_m\nabla u^s) = f^s \quad \mbox{in} \,\, \Omega_m.
\end{equation}

Subtracting (\ref{eqn:PE_sr}) from (\ref{eqn:PBE_splitting}) and using (\ref{eqn:Boltzmann_lr}), 
we obtain the nonlinear regularized PBE (NRPBE) as follows
\begin{equation}\label{eqn:RPBE_decomp}
-\nabla \cdot(\epsilon\nabla u^r(\bar{x})) + \bar{\kappa_2}^2(\bar{x})\sinh(u^r(\bar{x})) 
   = f^l,   \quad \mbox{in} \,\, \Omega,
\end{equation}
subject to the Dirichlet boundary conditions in (\ref{eq:DH_solution}). The total solution to the 
NRPBE is therefore, obtained by $u(\bar{x}) = u^s(\bar{x}) + u^r(\bar{x})$.

\section{Numerical approach to solving the NRPBE} \label{sec:RPBE_numerical_approach}
\subsection{Iterative solution of the NRPBE}
Let us consider a physical domain $\Omega \subset \mathbb{R}^3$ with boundary $\partial \Omega$, and 
a parameter domain $\mathcal{P} \subset \mathbb{R}$ which represents the variation in ionic strength 
$I = 1/2\sum_{j=1}^{N_{ions}}c_jz_j^2$, which is a function of the ionic concentration $c_i$, of the 
salt solution. It resides in $\bar{k}^2 = {8\pi e^2 I}/{1000\epsilon k_BT}$. One standard way of 
solving the NRPBE in (\ref{eqn:RPBE_decomp}) is that it is first discretized in space to obtain a nonlinear 
system in matrix-vector form 
\begin{equation}\label{eq:FOM}
 A(u_{\mathcal{N}}^r(\mu)) = b^r(\mu),
 \quad \mu \in \mathcal{P},
\end{equation}
where $A(u_{\mathcal{N}}^r(\mu)) \in \mathbb{R}^{\mathcal{N}\times \mathcal{N}}$, $b^r(\mu) 
\in \mathbb{R}^{\mathcal{N}}$, $\mu = I \in \mathcal{P}$, and $u_{\mathcal{N}}^r(\mu)$ is the discretized 
solution vector.

Then system (\ref{eq:FOM}) can be solved using several existing techniques. For 
example, nonlinear relaxation methods have been implemented in the Delphi software 
\cite{Rocchia_2001}, the nonlinear conjugate gradient (CG) method has been implemented in University 
of Houston Brownian Dynamics (UHBD) software \cite{Brock_1992}, the nonlinear multigrid (MG) method 
\cite{Oberoi_1993} and the inexact Newton method are available in the adaptive 
Poisson-Boltzmann solver (APBS) software \cite{Holst:95}.  

In this study, we apply a different approach of solving (\ref{eqn:RPBE_decomp})
\cite{Mirzadeh:13, Shestakov:2002, Ji:2018}. In particular, an iterative approach is first applied to the 
continuous NRPBE in (\ref{eqn:RPBE_decomp}), where at the $(n+1)$st iteration step, the NRPBE is approximated by 
a linear equation via the Taylor series truncation. The expansion point of the Taylor series is the 
continuous solution $(u^r(\mu))^n$ at the nth iteration step.

Consider $(u^r(\mu))^n$ as the approximate solution at the $n$th iterative step, then the nonlinear 
term $\sinh((u^r(\mu))^{n+1})$ at the $(n+1)$st step is approximated by its truncated Taylor series 
expansion as follows
\begin{equation}\label{eq:Taylor_expansion_sinh}
\sinh((u^r(\mu))^{n+1}) \approx \sinh((u^r(\mu))^n) 
+ ((u^r(\mu))^{n+1} - (u^r(\mu))^n)\cosh((u^r(\mu))^n).
\end{equation}
Substituting the approximation (\ref{eq:Taylor_expansion_sinh}) into (\ref{eqn:RPBE_decomp}), we 
obtain
\begin{multline}\label{eq:PBE_approx_sinh}
-\nabla\cdot(\epsilon(\bar{x})\nabla (u^r(\mu))^{n+1}) 
+ \bar{\kappa}^2(\bar{x})\cosh((u^r(\mu))^n)(u^r(\mu))^{n+1} 
= -\bar{\kappa}^2(\bar{x})\sinh((u^r(\mu))^n) \\
+ \bar{\kappa}^2(\bar{x})\cosh((u^r(\mu))^n)(u^r(\mu))^n + b^r(\mu).
\end{multline}
The equation in (\ref{eq:PBE_approx_sinh}) is linear, and can then be numerically solved by first applying 
spatial discretization. In this regard, we first define
\begin{equation}\label{eqn:hyperbolic_cosine_vec}
\cosh\odot u_{\mathcal{N}}^r(\mu) =: w = \begin{bmatrix}
		w_1 \\
		w_2 \\
		\vdots \\
		w_{\mathcal{N}}
		\end{bmatrix},
\end{equation}
where $\odot$ is the elementwise operation on a vector. 

Then, we construct the corresponding diagonal matrix from (\ref{eqn:hyperbolic_cosine_vec}) of the 
form \[B = \mbox{diag}(w_1, w_2, \ldots, w_{\mathcal{N}}).\] 
Finally, we obtain the following iterative linear system 
\begin{equation}\label{eqn:Nonaffine_form_iterative_FOM}
 A_1(u_{\mathcal{N}}^r(\mu))^{n+1} + \mu A_2B^{n}(u_{\mathcal{N}}^r(\mu))^{n+1} 
 = -\mu A_2\sinh\odot(u_{\mathcal{N}}^r(\mu))^{n} + \mu A_2B^{n}(u_{\mathcal{N}}^r(\mu))^{n} + b_1^r + b_2(\mu),
\end{equation}
where $A_1$ is the Laplacian matrix and $A_2$ is a diagonal matrix containing the net $\bar{\kappa}^2$ 
function (i.e., $\bar{\kappa}^2/\mu$). Note that the diagonal matrix $B^{n}$ changes at each iteration step, 
therefore, it cannot be precomputed. The vectors $b_1^r$ and $b_2(\mu)$ are the regularized approximation 
of the Dirac delta distributions and the nonaffine (in parameter) Dirichlet boundary conditions, respectively.

Let 
\begin{equation}\label{eqn:affine_A_iter}
 A(\cdot) = A_1 + \mu A_2B^n(\cdot)
\end{equation} 
and
\begin{equation}\label{eqn:affine_F_iter}
 F : \mbox{right-hand side of} \: (\ref{eqn:Nonaffine_form_iterative_FOM}),
\end{equation} 
we obtain
\begin{equation}\label{eq:PBE_system_iterative}
A((u_{\mathcal{N}}^r(\mu))^n)(u_{\mathcal{N}}^r(\mu))^{n+1} = F((u_{\mathcal{N}}^r(\mu))^n), 
\quad n = 0,1, \ldots.
\end{equation} 

Then, at each iteration, system (\ref{eq:PBE_system_iterative}) is a linear system w.r.t. 
$(u^r_{\mathcal{N}})^{n+1}$, which can be solved by any linear system solver of choice. In this study, we 
employ the aggregation-based algebraic multigrid method (AGMG) \footnote{\textbf{AGMG} implements an 
aggregation-based algebraic multigrid method, which solves algebraic systems of linear equations, and is 
expected to be efficient for large systems arising from the discretization of scalar second order elliptic 
PDEs \cite{Notay:2010}.} \cite{Notay:2010}. \Cref{alg:Iterative_NRPBE} summarizes the detailed iterative 
approach for solving (\ref{eq:PBE_system_iterative}). This approach of first linearization, then discretization 
is shown to be more efficient than the standard way of first discretization and then linearization, via, for 
example, the Newton iteration. The advantage of the proposed approach is that it avoids computing the Jacobian 
of a huge matrix. It is observed that it converges faster than the standard Newton approach. 

\begin{algorithm}[t]
  \caption{Iterative solver for the NRPBE}\label{alg:Iterative_NRPBE}
    \begin{algorithmic}[1]
    \Require Initialize the potential $(u_{\mathcal{N}}^r(\mu))^0$, e.g., $(u_{\mathcal{N}}^r(\mu))^0 = 0$ 
	      and the tolerance $\delta^0 = 1$. 
    \Ensure The converged NRPBE solution $(u_{\mathcal{N}}^r(\mu))^n$ at $\delta^n \leq \textrm{tol}$.
    \While{$\delta^n \geq \textrm{tol}$}
    \State Solve the linear system (\ref{eq:PBE_system_iterative}) for 
      $(u_{\mathcal{N}}^r(\mu))^{n+1}$ using AGMG.
    \State $\delta^{n+1} \gets \|(u_{\mathcal{N}}^r(\mu))^{n+1} - (u_{\mathcal{N}}^r(\mu))^n\|_2$.
    \State $(u_{\mathcal{N}}^r(\mu))^n \gets (u_{\mathcal{N}}^r(\mu))^{n+1}$.
    \EndWhile
    \State \textbf{end while}
  \end{algorithmic}
\end{algorithm}

\subsection{The reduced basis method for the NRPBE} \label{sec:RPBE_RBM}

The Reduced Basis Method (RBM) is an example of popular projection-based parametrized model order 
reduction (PMOR) techniques in which the parameter dependence of the PDE solution is exploited by 
snapshots (high-fidelity solutions) determined over the parameter domain \cite{morHest16}. Their core 
objective is to construct a parametric reduced order model (ROM) of low dimension, which accurately 
approximates the original full order model (FOM) or high-fidelity model of high dimension over varying 
parameter values \cite{morBenGW15, morEft11,morRozHP08}. Other PMOR techniques include proper 
orthogonal decomposition (POD) \cite{morVol13} and multi-moment matching techniques \cite{morBenF14}, among 
others \cite{morBenGW15}. 

The RBM leverages an offline/online procedure to ensure an accurate approximation of the high-fidelity 
solution at extremely low computational costs. It is widely applicable in real-time context such as 
sensitivity analysis, multi-model simulation, as well as many-query scenarios, e.g., uncertainty 
quantification and optimal control. For a thorough review of PMOR techniques, see \cite{morBenGW15}. 

It is prohibitively expensive to solve the $\mathcal{N}\times \mathcal{N}$ system in 
(\ref{eq:PBE_system_iterative}) for an accurate approximation of $u_{\mathcal{N}}^r(\mu)$ because 
the dimension $\mathcal{N}$ ranges between $\mathcal{O}(10^6)$ and $\mathcal{O}(10^8)$ for typical 
biomolecules of interest. Therefore, we exploit the RBM to significantly reduce the computational 
costs by projecting (via Galerkin) the FOM (\ref{eq:PBE_system_iterative}) onto a low dimensional subspace 
(the reduced basis space) which preserves the parametric properties and important system configurations 
of the FOM. The resultant ROM, of greatly reduced dimension $N \ll \mathcal{N}$, provides an accurate 
surrogate approximation of $(u_{\mathcal{N}}^r(\mu))^n$, $n = 0,1, \ldots$. 

RBM is based on the assumption that the solution manifold $\mathcal{M}^{\mathcal{N}}$
\begin{equation}\label{Manifold}
 \mathcal{M}^{\mathcal{N}} = \{ u_{\mathcal{N}}^r(\mu) : \mu \in \mathcal{P}\},
\end{equation}
is of low dimension. The reduced basis space, which is the space spanned by the snapshots 
$u_{\mathcal{N}}^r(\mu)$, corresponding to a set of parameter samples, is hierarchically constructed 
from the solution manifold (\ref{Manifold}), using the greedy procedure summarized in  
\Cref{alg:greedy}. The RB space
\begin{equation}\label{RB_space}
 \mbox{range}(V) = \mbox{span}\{ u_{\mathcal{N}}^r(\mu_1), \ldots, u_{\mathcal{N}}^r(\mu_l)\}, 
 \quad  \mu_1, \ldots, \mu_l \in \mathcal{P},
\end{equation}
is nested (hierarchical) in the sense that the previous basis set is a subset of the next until 
convergence, i.e.,
 \begin{equation}\label{RB_space_nested}
 \mbox{range}(V_1) \subset \mbox{range}(V_2) \subset \cdots \subset \mbox{range}(V_N),
\end{equation}
where $N \ll \mathcal{N}$ is the dimension of the ROM.

The residual in \Cref{alg:greedy} is derived from (\ref{eq:PBE_system_iterative}) and the ROM 
solution $(\widehat{u}_{\mathcal{N}}^r(\mu))^n = V_N(u_N^r(\mu))^n$ lifted into the high-fidelity space 
of dimension $\mathcal{N}$, i.e.,
\begin{equation}\label{eqn:residual}
 r_N((\widehat{u}_{\mathcal{N}}^r(\mu))^{n+1}) = F((\widehat{u}_{\mathcal{N}}^r(\mu))^n) 
		- A((\widehat{u}^r(\mu))^n)(\widehat{u}_{\mathcal{N}}^r(\mu))^{n+1}. 
\end{equation}

\begin{algorithm}[t]
    \caption{Greedy algorithm}\label{alg:greedy}
      \begin{algorithmic}[1]
	\Require Training set $\Xi:= \{\mu_1, \ldots, \mu_l\} \subset \mathcal{P}$, tolerance 
	$\epsilon_0 = 1$, and potential $(u_N^r(\mu))^0$.
	\Ensure RB basis represented by $V$ and the ROM in (\ref{eq:ROM}).
	\State Choose $\mu^* \in \Xi$ arbitrarily.
	\State Solve (\ref{eqn:RPBE_decomp}) for $u_{\mathcal{N}}^r(\mu^*)$ using 
		\Cref{alg:Iterative_NRPBE}.
	\State $V_1 = [u_{\mathcal{N}}^r(\mu^*)]$, $N = 1$.
	\State Orthonormalize $V_1$.
	\While{$\max\limits_{\mu\in \Xi}\Delta_N(\mu) \geq \epsilon$}
	\State Compute $u_N^r(\mu)$ from (\ref{eq:ROM}) using \Cref{alg:Iterative_RROM}, and calculate 
		$\Delta_N(\mu) = \Arrowvert r_N(\widehat{u}_{\mathcal{N}}^r(\mu))\Arrowvert_2$ 
	\StatexIndent[1] in (\ref{eqn:residual}), $\forall \, \mu \in \Xi$.
	\State $\mu^* = \mbox{arg}\max\limits_{\mu\in \Xi}\Delta_N(\mu)$.
	\State Solve (\ref{eqn:RPBE_decomp}) for $u_{\mathcal{N}}^r(\mu^*)$.
	\State $V_{N+1} \gets [V_N \quad u_{\mathcal{N}}^r(\mu^*)]$.
	\State Orthonormalize the columns of $V_{N+1}$.
	\State $N \gets N + 1$.
	\EndWhile
	\State \textbf{end while}
      \end{algorithmic}
  \end{algorithm}

The ROM for the system (\ref{eq:PBE_system_iterative}), is therefore, formulated 
as follows. Given any $\mu \in \mathcal{P}$, and an initial guess $(u_N^r(\mu))^0 \in \mathbb{R}^N$, 
the RB approximation $(u_N^r(\mu))^{n+1}$, at the future iteration step $n+1$ satisfies the equation
\begin{equation}\label{eq:ROM}
 A_N((u_N^r(\mu))^{n})(u_N^r(\mu))^{n+1} = F_N((u_N^r(\mu))^{n}), \quad n = 0,1, \ldots, 
\end{equation}
where $(u_N^r(\mu))^{0}$ is the zero initial guess in this study and $A_N$ and $F_N$ are defined 
explicitly as 
\[A_N := \widehat{A}_1(u_N^r(\mu))^{n+1} 
+ \mu \widehat{A}_2(\widetilde{B}V_N)^n(u_N^r(\mu))^{n+1},\] and
\[F_N := -\mu \widehat{A}_2\sinh\odot(\widehat{u}_{\mathcal{N}}^r(\mu))^n 
 + \mu \widehat{A}_2(\widetilde{B}V_N)^n(u_N^r(\mu))^{n} + b_N^r + V_N^Tb_2(\mu),\]
where $\widetilde{B} = \mbox{diag}(\widetilde{w}_1, \widetilde{w}_2, \ldots, \widetilde{w}_{\mathcal{N}})$ 
and 
\begin{equation}\label{eqn:hyperbolic_cosine_vec2}
\cosh\odot \widehat{u}_{\mathcal{N}}^r(\mu) =: \widetilde{w} = \begin{bmatrix}
		\widetilde{w}_1 \\
		\widetilde{w}_2 \\
		\vdots \\
		\widetilde{w}_{\mathcal{N}}
		\end{bmatrix}.
\end{equation}

The resuting ROM is given by
\begin{multline}\label{eq:ROM_detail}
\widehat{A}_1(u_N^r(\mu))^{n+1} 
+ \mu \widehat{A}_2(\widetilde{B}V_N)^n(u_N^r(\mu))^{n+1} 
 = -\mu \widehat{A}_2\sinh\odot(\widehat{u}_{\mathcal{N}}^r(\mu))^n \\
 + \mu \widehat{A}_2(\widetilde{B}V_N)^n(u_N^r(\mu))^{n} + b_N^r + V_N^Tb_2(\mu),
\end{multline}
where $(u_N^r(\mu))^{n+1}$ is the unknown solution to the ROM.

The reduced matrices $\widehat{A}_1 := V_N^T A_1 V_N$ and $\widehat{A}_2 := V_N^T A_2$ and the reduced vector 
$b_N^r = V_N^Tb_1^r$ (see (\ref{eqn:Nonaffine_form_iterative_FOM})) are determined via projection with the 
RB basis $V_N$ and can be precomputed in the offline phase of the greedy algorithm. However, the matrix 
$(\widetilde{B}V_N)^n$ and vector $b_2(\mu)$ are updated and or changed at each iteration and for varying parameter 
values, respectively, hence, they cannot be precomputed. This leads to a partial offline-online decomposition 
scenario, whereby Galerkin projections to some terms have to be computed in the online phase. 

Note that $V_N^Tb_2(\mu)$ in (\ref{eq:ROM_detail}) is computed by first evaluating a long vector 
$b_2(\mu)$, then projecting it onto the low dimensional space $N$ using $V_N$. This is time consuming 
when $b_2(\mu)$ needs to be evaluated many times for many values of $\mu$. In \Cref{sec:RPBE_RROM}, 
we propose to apply DEIM to further reduce the computational complexity of $V_N^Tb_2^r(\mu)$. 
Details can be found in \cite{KwBFM2017}, where DEIM was applied to PBE problem.

\subsubsection{Computational complexity of the regularized reduced order model} \label{sec:RPBE_RROM}

It is well known that another key assumption of the RBM, besides the low dimensionality of the solution 
manifold, is the parameter affine property, which ensures the efficiency of the offline-online 
decomposition by eliminating the dependency of the ROM on the dimension $\mathcal{N}$ of the truth 
high-fidelity FOM \cite{morHest16}. However, note that on the one hand, 
(\ref{eqn:Nonaffine_form_iterative_FOM}) is actually parameter nonaffine with respect to the 
Yukawa-type boundary conditions, represented by $F$ in (\ref{eqn:affine_F_iter}). On the other hand, 
the matrix $A_2$ requires updates at each iteration, hence Galerkin projections are unavoidable in the 
online phase. 

In this study, we apply DEIM to the parametric nonaffine boundary conditions, the term $b_2(\mu)$. The main 
idea of DEIM is to significantly reduce the computational complexity of the nonaffine function 
by interpolation, whereby only a few entries are computed \cite{morChaS10}. Before invoking DEIM, 
snapshots of the nonaffine function $b_2(\mu)$ must be computed for a set of parameter $\mu$ 
in the training set $\Xi=\{\mu_1, \ldots, \mu_l\}\subset \mathcal{P}$ and the snapshot matrix,
\begin{equation}
 G=[b_2(\mu_1), \ldots, b_2(\mu_l)] \in \mathbb{R}^{\mathcal{N}\times l},
\end{equation}
is constructed.

\begin{figure}[t]
  \centering
  \captionsetup{width=\linewidth}
    \includegraphics[height=8.0cm]{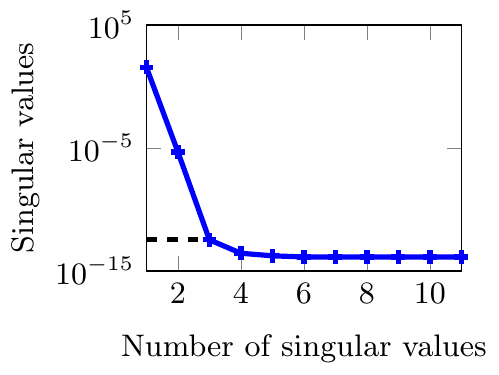}
  \caption{Decay of singular values of $\Sigma$ in (\ref{eq:svd}).}
  \label{fig:Singular_values_1FAS}
\end{figure}

Then, the singular value decomposition (SVD) of $G$ is computed, 
\begin{equation}\label{eq:svd}
 G = U_G\Sigma W^T,
\end{equation}
where $U_G \in \mathbb{R}^{\mathcal{N}\times l}$, $\Sigma \in \mathbb{R}^{l\times l}$, and 
$W \in \mathbb{R}^{l\times l}$. Note that the matrices $U_G$ and $W$ are orthogonal, i.e., 
$(U_G)^TU_G=W^TW = I_l$, $I_l \in \mathbb{R}^{l \times l}$ and $\Sigma = \textrm{diag}(\sigma_1, 
\ldots, \sigma_l)$, with $\sigma_1 \geq \ldots \geq \sigma_l \geq 0$. \Cref{fig:Singular_values_1FAS} 
shows the decay of singular values of $\Sigma$ for the protein \textit{fasciculin 1}. We truncate the 
singular values being smaller than $10^{-13}$ and retain $r=3$ singular values.   

\setlength{\fwidth}{6.0cm}
\setlength{\fheight}{4.5cm}
 
DEIM seeks to approximate the function $b_2(\mu)$ with the linear combination of the basis vectors 
$U_G=[u_1^G, \ldots, u_r^G] \in \mathbb{R}^{\mathcal{N} \times r}$, i.e.
\begin{equation}
 b_2(\mu) \approx U_Gc(\mu),
\end{equation}
where $c(\mu)\in \mathbb{R}^r$ is the corresponding coefficient vector, which is determined under the 
assumption that $U_Gc(\mu)$ interpolates $b_2(\mu)$ at $r$ selected interpolation points, then,
\begin{equation}\label{eq:Overdet_sys}
 P^Tb_2(\mu) = P^TU_Gc(\mu),
\end{equation}
where $P$ is an index matrix given by
\begin{equation}
 P = [e_{\wp_1}, \ldots, e_{\wp_r}] \in \mathbb{R}^{\mathcal{N} \times r},
\end{equation} 
which consists of unit vectors $e_{\wp_i}$, $i = 1, \ldots, r$. Here, the indices $\wp_i$, are the 
DEIM interpolation points which are selected iteratively with the greedy iteration as presented in 
\Cref{alg:DEIM}. 

Suppose that $P^TU_G \in \mathbb{R}^{r \times r}$ is nonsingular, then $c(\mu)$ can be determined 
from the overdetermined system (\ref{eq:Overdet_sys}) by
\begin{equation}
c(\mu) = (P^TU_G)^{-1}P^Tb_2(\mu).
\end{equation}

\begin{algorithm}[t]
  \caption{DEIM algorithm \cite{morChaS10, Feng2016}}\label{alg:DEIM}
  \begin{algorithmic}[1]
      \Require POD basis $\{u_i^G\}_{i=1}^r$ for $G$ in equation (\ref{eq:svd})). 
      \Ensure DEIM basis $U_G$ and indices $\vec{\wp} = [\wp_1,\ldots,\wp_r]^T \in \mathbb{R}^r$.
	\State $\wp_1$ = $\textrm{arg} \max\limits_{j \in \{1, \ldots, \mathcal{N}\}}\lvert u_{1j}^G 
		\rvert$, where $u_1^G = (u_{11}^G, \ldots, u_{1\mathcal{N}}^G)^T$. 
	\State $U_G = [u_1^G]$, $P = [e_{\wp_1}], \vec{\wp} = [\wp_1]$.
	  \For{\texttt{i = 2 to r}}
	  \State Solve $(P^TU_G)\alpha = P^Tu_i^G$ for $\alpha$, where 
	  $\alpha = (\alpha_1, \ldots, \alpha_{i-1})^T$.
	  \State $r_i = u_i^G - U_G \alpha$.
	  \State $\wp_i$ = $\textrm{arg} \max\limits_{j \in \{1, \ldots, 
		  \mathcal{N}\}}\lvert r_{ij} \rvert$, where $r_i = (r_{i1}, 
		  \ldots, r_{i\mathcal{N}})^T$. 
	  \State $U_G \gets [U_G \ u_i^G]$, $P \gets [P \ e_{\wp_i}]$, $\vec{\wp} \gets 
		\begin{bmatrix}
		\vec{\wp} \\
		\wp_i
		\end{bmatrix}$.
	  \EndFor
	\State \textbf{end for}
  \end{algorithmic}
\end{algorithm}

Therefore, the function $b_2(\mu)$ in (14) can be approximated as
\begin{equation}\label{eq:DEIM_appro}
 b_2(\mu) \approx U_Gc(\mu) = U_G(P^TU_G)^{-1}P^Tb_2(\mu).
\end{equation}
 
The ROM in (\ref{eq:ROM_detail}) with DEIM approximation becomes
\begin{multline}\label{eq:ROM_detail_DEIM}
 \widehat{A}_1(u_N^r(\mu))^{n+1} + \mu \widehat{A}_2(\widetilde{B}V_N)^n(u_N^r(\mu))^{n+1} 
 = -\mu \widehat{A}_2\sinh\odot(\widehat{u}_{\mathcal{N}}^r(\mu))^n 
 + \mu \widehat{A}_2(\widetilde{B}V_N)^n(u_N^r(\mu))^{n} \\+ b_N^r(\mu) 
 + V_N^TU_G(P^TU_G)^{-1}P^Tb_2(\mu).
\end{multline}

The interpolant $V_N^TU_G(P^TU_G)^{-1}P^Tb_2(\mu)$ can be computed a lot cheaper than $V_N^Tb_2(\mu)$ 
because $V_N^TU_G(P^TU_G)^{-1}$ can be precomputed independently of the parameter $\mu$. Only 
those entries in $b_2(\mu)$ that correspond to the interpolation indices $\wp_i, i = 1, \ldots, 
r$, $r \ll \mathcal{N}$, i.e., $P^T b_2(\mu)$ can be computed instead of the entire $\mathcal{N}$ 
entries in $b_2(\mu)$. This saves significant computational efforts when $b_2(\mu)$ needs to be repeatedly 
computed for different values of $\mu$.

Note that at each iteration, only a small ROM in (\ref{eq:ROM_detail_DEIM}) is solved. With its small 
size $N \ll \mathcal{N}$, the system (\ref{eq:ROM_detail_DEIM}) can be solved using a direct solver 
rather than the iterative solver (AGMG), which is applied to to the FOM in (\ref{eq:PBE_approx_sinh}). 
The iterative approach of obtaining an approximate solution $V_N (u_N^r(\mu))^{n+1}$ to 
(\ref{eq:FOM}) using the ROM (\ref{eq:ROM_detail_DEIM}) is summarized in \Cref{alg:Iterative_RROM}.

\begin{algorithm}[t]
  \caption{Iterative solver for the regularized ROM in (\ref{eq:ROM_detail_DEIM})}\label{alg:Iterative_RROM} 
    \begin{algorithmic}[1]
    \Require Initialize the potential $(u_N^r(\mu))^0$, e.g., $(u_N^r(\mu))^0 = 0$, 
	       tolerance tol $>0$, and $\delta^0 = 1$. 
    \Ensure The converged ROM solution $(u_N^r(\mu))^n$ at $\delta^n \leq \textit{tol}$.
    \State Precompute $\widehat{A}_1$, $b_N^r$ in (\ref{eq:ROM_detail_DEIM}) and $U_G$ and $\vec{\wp}$ 
	    in \Cref{alg:DEIM}.
    \While{$\delta^n \geq \textit{tol}$}
    \State Assemble the ROM in (\ref{eq:ROM_detail_DEIM}) using the precomputed quantities in Step 1.
    \State Solve the regularized ROM (\ref{eq:ROM_detail_DEIM}) for $(u_N^r(\mu))^{n+1}$.
    \State $\delta^{n+1} \gets \|(u_N^r(\mu))^{n+1} - (u_N^r(\mu))^n\|_2$.
    \State $(u_N^r(\mu))^n \gets (u_N^r(\mu))^{n+1}$.
    \EndWhile
    \State \textbf{end while}
  \end{algorithmic}
\end{algorithm}

\begin{remark}
 The total electrostatic potential is obtained by lifting the reduced order long-range surrogate 
 solution into the high-fidelity space $\mathcal{N}$ and adding to the parameter independent 
 analytically precomputed short-range component ${\bf P}_s$ in (\ref{eq:CCT_tensor}), i.e.,
 \begin{equation}
  u(\mu) = {\bf P}_s + \widehat{u}_{\mathcal{N}}^r(\mu),
 \end{equation}
  where $\widehat{u}_{\mathcal{N}}^r(\mu) = V_Nu_N^r(\mu)$.
\end{remark}

\subsection{The reduced basis method for the classical NPBE} \label{sec:NPBE_RBM}

In this section, we apply RBM to the classical nonlinear PBE (NPBE), and compare this version 
in \Cref{sec:Numer_Tests} with the suggested approach from \Cref{sec:RPBE_RBM}. We begin by considering the FOM of the classical NPBE in (\ref{eqn:PBE}) 
after discretization in space, i.e.,
\begin{equation}\label{eq:FOM_classic}
 A(u^{\mathcal{N}}(\mu)) = f(\mu), \quad \mu \in \mathcal{P},
\end{equation}
where $f(\mu)$ includes both the singular sources from the right-hand side of (\ref{eqn:PBE}) and the 
parameter non-affine Dirichlet boundary conditions from (\ref{eq:DH_solution}). The corresponding 
classical ROM is defined as 
\begin{equation}\label{eqn:ROM_class1}
 \widehat{A}(u_N(\mu)) = \widehat{f} (\mu),
\end{equation}
where $\widehat{A}(u_N(\mu)) = V_N^T A(V_Nu_N(\mu))$ and $\widehat{f} = V_N^T f$. Here, $V_N$ can be constructed 
using the greedy algorithm in \Cref{alg:greedy} by replacing the snapshots in Step 3 and Step 9 with the 
solutions to (\ref{eq:FOM_classic}). 

Note that the FOM of the classical NPBE is solved iteratively in a similar way like the NRPBE using 
\Cref{alg:Iterative_NRPBE}. The corresponding iterative form of (\ref{eq:FOM_classic}) is given by
\begin{equation}\label{eq:FOM_classic_iterative}
 A_1(u^{\mathcal{N}}(\mu))^{n+1} + \mu A_2B_2^{n}(u^{\mathcal{N}}(\mu))^{n+1} 
 = -\mu A_2\sinh\odot(u^{\mathcal{N}}(\mu))^{n} + \mu A_2B_2^{n}(u^{\mathcal{N}}(\mu))^{n} + f(\mu),
\end{equation}
where all the quantities except $B_2$, $f$ and the solution $(u^{\mathcal{N}}(\mu))^{n+1}$ are equivalent to 
those in the NRPBE (\ref{eqn:Nonaffine_form_iterative_FOM}). Here $B_2$ is defined as
\[B_2 = \mbox{diag}(v_1, v_2, \ldots, v_{\mathcal{N}}),\]
and is constructed from 
\begin{equation}\label{eqn:hyperbolic_cosine_vec3}
\cosh\odot u^{\mathcal{N}}(\mu) =: v = \begin{bmatrix}
		v_1 \\
		v_2 \\
		\vdots \\
		v_{\mathcal{N}}
		\end{bmatrix}.
\end{equation}

The ROM of (\ref{eq:FOM_classic_iterative}) is straightforward, i.e., given any $\mu \in 
\mathcal{P}$, and an initial potential distribution $(u_N(\mu))^0$, the RB approximation 
$(u_N(\mu))^{n+1}$, at the subsequent iteration steps $n+1$ satisfies 
\begin{multline}\label{eq:ROM_classical1}
\widehat{A}_1(u_N(\mu))^{n+1} + \mu \widehat{A}_2(\widetilde{B_2}V_N)^n(u_N(\mu))^{n+1} 
 = -\mu \widehat{A}_2\sinh\odot(\widehat{u}^{\mathcal{N}}(\mu))^n \\
 + \mu \widehat{A}_2(\widetilde{B_2}V_N)^n(u_N(\mu))^{n} + \widehat{f}(\mu),
\end{multline}
where $(\widehat{u}^{\mathcal{N}}(\mu))^n = V_N(u_N(\mu))^n$ and
$\widetilde{B_2} = \mbox{diag}(\widetilde{v}_1, \widetilde{v}_2, \ldots, \widetilde{v}_{\mathcal{N}})$ is 
constructed from
\begin{equation}\label{eqn:hyperbolic_cosine_vec4}
\cosh\odot \widehat{u}^{\mathcal{N}}(\mu) =: \widetilde{v} = \begin{bmatrix}
		\widetilde{v}_1 \\
		\widetilde{v}_2 \\
		\vdots \\
		\widetilde{v}_{\mathcal{N}}
		\end{bmatrix}.
\end{equation}
The process of iteratively solving (\ref{eq:ROM_classical1}) is similar to that of (\ref{eq:ROM}), 
which is provided in \Cref{alg:Iterative_RROM}.

\section{Numerical results}
\label{sec:Numer_Tests}

Consider $n^{\otimes 3}$ 3D uniform Cartesian grids, in a cubic domain $[a,b]^3$, for computing the 
reduced basis approximation of the NRPBE on a modest PC which possesses the following specifications: 
Intel (R) Core (TM) $i7-4790$ CPU @ 3.60GHz with 8GB RAM. In this study, the NRPBE is discretized by 
the finite difference method (FDM) to obtain the FOM and the numerical computations are implemented 
in MATLAB, version R2017b.

In the numerical tests, the molecular charge density function (singular source term) for the classical NPBE 
and the regularized Dirac density function for the NRPBE are obtained from PQR 
\footnote{A PQR (or Position, charge (Q), and Radius) file is a protein data bank (PDB) 
file with the temperature and occupancy columns replaced by columns containing the per-atom 
charge (Q) and radius (R) using the pdb2pqr software. PQR files are used in several computational 
biology packages, including APBS \cite{Bakersept2001}.} files which are generated from the following 
biomolecules with varying sizes that depend on the number of atoms:
\begin{enumerate}[label=(\alph*)]
\item The \textit{acetazolamide} molecule consisting of $18$ atoms, which is used as a ligand in the human 
carbonic anhydrase (hca) protein-ligand complex for the calculation of the binding energy 
\cite{Holst:93,Vergara-Perez2016}.
 \item \textit{fasciculin 1}, an anti-acetylcholinesterase toxin from green mamba snake venom 
 \cite{DuMaBoFo:92} consisting of 1228 atoms.
\item A 180-residue cytokine solution NMR structure of a murine-human chimera 
of leukemia inhibitory factor (LIF) \cite{HinMauZhaNic:98} consisting of 2809 atoms.
\end{enumerate}

\begin{remark}
 Since the solution of the PBE has a slow polynomial decay in $1/\|\bar{x}\|$, it is paramount that 
 large domains, approximately $3$-times the size of the biomolecule be used in order to accurately 
 approximate the boundary conditions \cite{Holst94}. In this regard, we use domains of lengths 
 $(32\mbox{\AA})^3$, $(60\mbox{\AA})^3$ and $(65\mbox{\AA})^3$, respectively, for the aforementioned 
 biomolecules. Here, $\mbox{\AA}$ denotes the angstrom unit of length.
\end{remark}

\begin{figure}[t]
\centering
\captionsetup{width=\linewidth}
\resizebox{17cm}{9.0cm}{%
\begin{tikzpicture}

\node (Short_range_protein) at (0,3.0) {\begin{tabular}{l}
  {\includegraphics[width=0.48\textwidth]{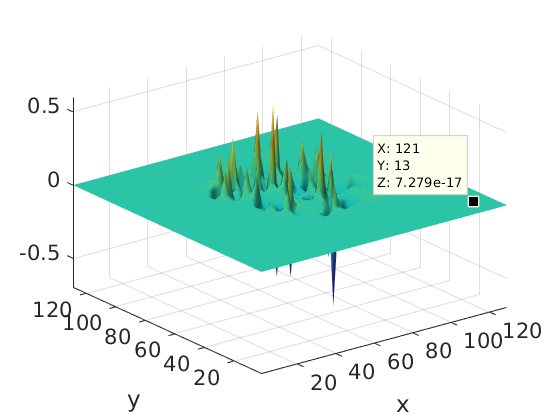}}
 \end{tabular}
 };
 
 \node (Long_range_protein) at (0,-3.0) {\begin{tabular}{l}
  {\includegraphics[width=0.48\textwidth]{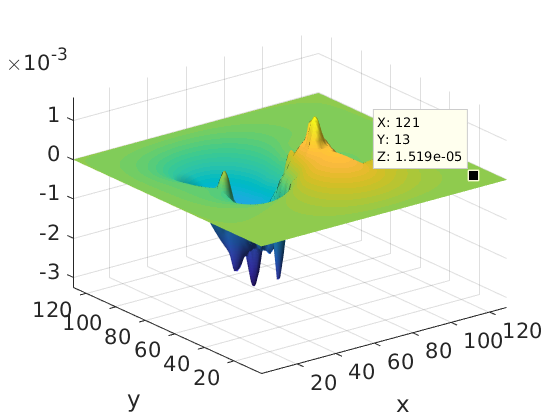}}
 \end{tabular}
 };
 
 \node (Full_pot_protein) at (8.5,0) {\begin{tabular}{l}
  {\includegraphics[width=0.48\textwidth]{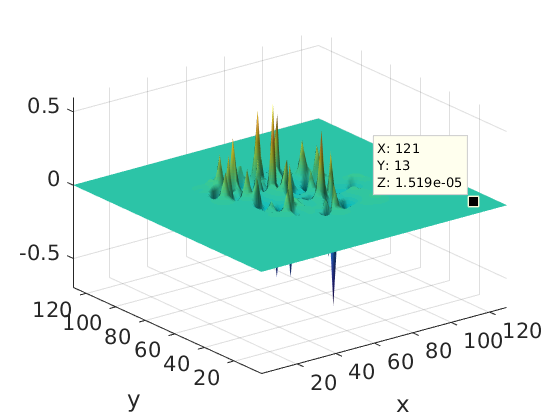}}
 \end{tabular}
 };
 \draw[->,draw=blue,thick] (3.5,0) -- (Full_pot_protein.west);
 
 \end{tikzpicture}}
 \caption{\label{fig:Protein_pot_RS}
 The long-range (bottom left), short-range (top left), and total electrostatic potentials (right) 
 for \textit{fasciculin 1}.}
\end{figure}

To begin with, we demonstrate the solution components of the full order model (FOM) of the NRPBE via the RS 
tensor format for the protein \textit{fasciculin 1} in case (b), in a uniform Cartesian grid of 
$129^{\otimes3}$ and a $60\mbox{\AA}$ domain length. \Cref{fig:Protein_pot_RS} shows the short- and 
long-range components of the target electrostatic potential, which are computed analytically from the CCT 
tensor (\ref{eq:CCT_tensor}), and numerically via the NRPBE in (\ref{eqn:RPBE_decomp}), respectively, and 
the corresponding total electrostatic potential. 

 The behaviour in the scaling among the three electrostatic potentials in \Cref{fig:Protein_pot_RS} shows 
 that the total potential on the right-hand side inherits the largest potential value of $0.5$ from the 
 short-range component (top left), while the smallest value of $\mathcal{O}(10^{-5})$ is obtained from the 
 long-range component (bottom left) towards the boundary of the domain.

\begin{remark}
 The main computational advantage of applying the RBM technique to the NRPBE is that the RB approximation 
 is only applied to the smooth long-range component of the potential, see \Cref{fig:Protein_pot_RS} 
 (bottom left), thereby avoiding the singularities inherent in the short-range component that are known to 
 cause numerical difficulties. Hence the resultant RB approximation is expected to be of higher accuracy.
\end{remark}

\subsection{Reduced basis results for the NRPBE}
\label{ssec:NRPBE_RBM} 

Here, we determine the accuracy and computational efficacy of the RBM approximating the high-fidelity 
solution to the NRPBE for biomolecular modeling. We set the solute and solvent dielectric coefficients 
as $\epsilon_m = 2$, and $\epsilon_s = 78.54$, respectively, and employ the parameter values from the 
training set $\Xi \in \mathcal{P} = [0.05,0.15]$ with a sample size of $l=11$, greedy tolerance 
$\mbox{tol}= 10^{-10}$ for \Cref{alg:greedy} to generate the projection matrix $V_N$. Furthermore the 
residual in (\ref{eqn:residual}) is used as an error estimator for the ROM in the greedy algorithm in 
\Cref{alg:greedy}.

First, we consider the NRPBE system generated by all the three cases (a), (b) and (c), in uniform Cartesian 
grids of $97^{\otimes 3}$ for case (a) and $129^{\otimes 3}$ for cases (b) and (c), respectively. We show in 
\Cref{table:Max_vs_true_error}, the decay of the maximal error estimator, defined as
\[\Delta_N^{\max}(\mu) = \max\limits_{\mu \in \mathcal{P}} \Arrowvert r_N(\widehat{u}_N^r;\mu)\Arrowvert_2,\] 
and the true error $\Arrowvert u_{\mathcal{N}}^r(\mu)-\widehat{u}_N^r(\mu)\Arrowvert_2$, during the 
greedy algorithm at the current RB dimension $i = 1,\dots, N$ for all of these cases. 

Note that the ROM provides highly accurate approximations, close to machine precision 
($\mathcal{O}(10^{-15})$) for the NRPBE as demonstrated by the true error in the second iteration. This is due 
to the smoothness of the long-range electrostatic potential, which enhances rapid and accurate model reduction 
process and facilitates, in general, low-rank approximation.  

\begin{table}[t]
\centering
\captionsetup{width=\linewidth}
\begin{tabular}{|c|c|c|c|c|c|}
  \hline
  \multirow{2}{*}{Biomolecule} & \multicolumn{2}{c|}{Error at iteration 1} 
  & \multicolumn{2}{c|}{Error at iteration 2} & \multirow{2}{*}{ROM dimension $N$}\\
    \cline{2-5}
   & $\Delta_N^{\max}(\mu)$ & True error & $\Delta_N^{\max}(\mu)$ & True error & \\ \hline
  Case (a) & 5.0573e-06 & 1.2719e-08 & 3.0339e-12 & 3.0395e-15 & 2 \\ \hline
  Case (b) & 1.0685e-05 & 8.9228e-08 & 3.6895e-12 & 2.0232e-14 & 2 \\ \hline
  Case (c) & 3.2610e-05 & 1.4510e-07 & 2.0573e-11 & 3.2015e-14 & 2 \\ \hline
  \end{tabular}
\caption{The comparison between the maximal error estimator $\Delta_N^{\max}(\mu)$ and the true error for 
the NRPBE during the greedy iteration at the current RB dimension $i = 1,\dots, N$ for the biomolecules in 
cases (a) to (c).}\label{table:Max_vs_true_error}
\end{table}

\begin{figure}[t]
\captionsetup{width=\linewidth}
\begin{center}
    \begin{tikzpicture}
    \begin{customlegend}[legend columns=-1, legend style={/tikz/every even column/.append 
    style={column sep=2.0cm}} , legend entries={True error, $\Delta_N(\mu)$, tol}]
    \addlegendimage{blue,solid, line width = 2pt}
    \addlegendimage{black,solid, line width = 2pt}
    \addlegendimage{black,dashed, line width = 2pt}
    \end{customlegend}
    \end{tikzpicture}
    \end{center}

    \begin{subfigure}[b]{0.3\textwidth}
      \centering
      \includegraphics[height=4.0cm]{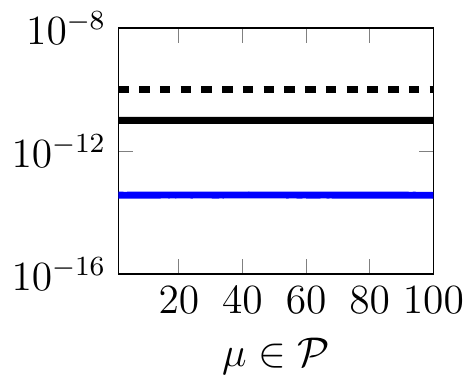}
      \caption{Error for case (a).}
      \label{fig:True_apost_compare_acet_current_rs}
    \end{subfigure}
    \hfill
    \begin{subfigure}[b]{0.3\textwidth}
        \centering
        \includegraphics[height=4.0cm]{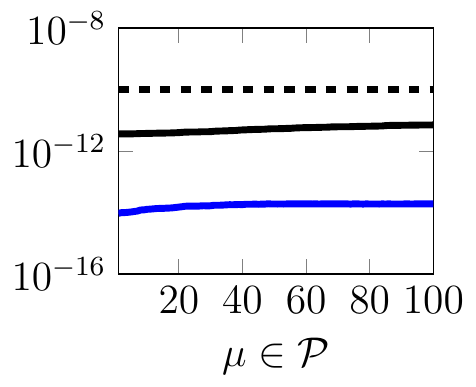}
        \caption{Error for case (b).}
        \label{fig:True_apost_compare_2809_current_rs}
    \end{subfigure}
    \hfill
    \begin{subfigure}[b]{0.3\textwidth}
        \centering
        \includegraphics[height=4.0cm]{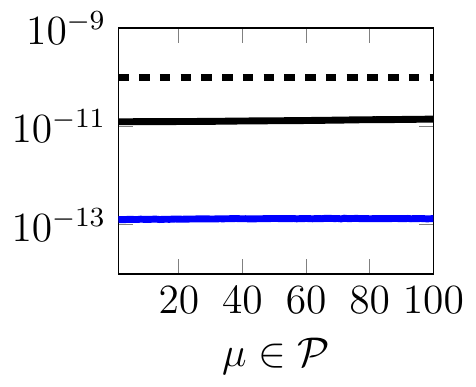}
        \caption{Error for case (c).}
        \label{fig:True_apost_compare_1FAS_final_rs}
    \end{subfigure}
    \caption{\label{fig:True_max_rand_compare_final_NRPBE}%
  Comparison between the error estimator and the true error for the NRPBE for the cases (a) to (c) 
  for the final ROM at $100$ random (varying) parameter values $\mu \in \mathcal{P}$.}
\end{figure}

Next, we validate the final ROM at $100$ random $\mu \in \mathcal{P}$ in 
\Cref{fig:True_max_rand_compare_final_NRPBE}. It is clear that the true error of the ROM is still below the 
tolerance for all $100$ $\mu \in \mathcal{P}$.

\setlength{\fwidth}{3.8cm}
\setlength{\fheight}{3.5cm}

\subsection{Comparison of the RB approximation accuracy between the NRPBE and the NPBE}
\label{ssec:NRPBE_vs_NPBE_RBM} 

\begin{figure}[t]
\captionsetup{width=\linewidth}
\begin{center}
    \begin{tikzpicture}
    \begin{customlegend}[legend columns=-1, legend style={/tikz/every even column/.append 
    style={column sep=2.0cm}} , legend entries={True error, $\Delta_N(\mu)$, tol}]
    \addlegendimage{blue,solid, line width = 2pt}
    \addlegendimage{black,solid, line width = 2pt}
    \addlegendimage{black,dashed, line width = 2pt}
    \end{customlegend}
    \end{tikzpicture}
    \end{center}

    \begin{subfigure}[b]{0.3\textwidth}
      \centering
      \includegraphics[height=4.0cm]{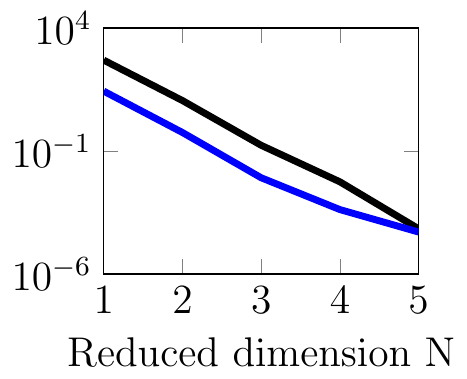}
      \caption{Error for case (a).}
      \label{fig:True_apost_compare_acet_current}
    \end{subfigure}
    \hfill
    \begin{subfigure}[b]{0.3\textwidth}
        \centering
        \includegraphics[height=4.0cm]{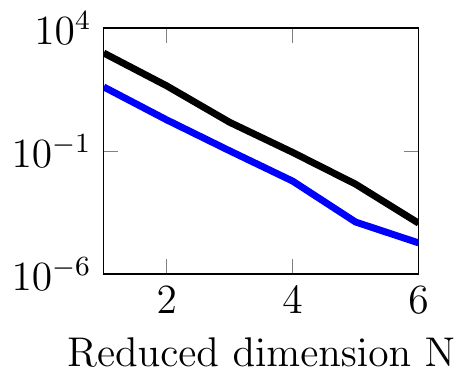} 
        \caption{Error for case (b).}
        \label{fig:True_apost_compare_1FAS_current}
    \end{subfigure}
      \hfill
    \begin{subfigure}[b]{0.3\textwidth}
        \centering
        \includegraphics[height=4.0cm]{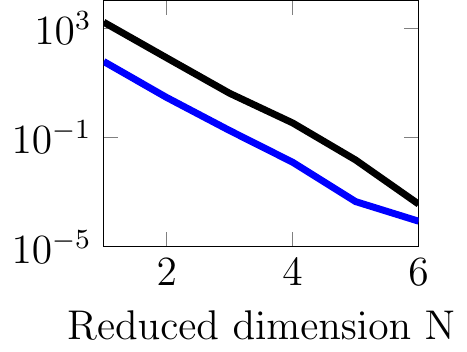}
        \caption{Error for case (c).}
        \label{fig:True_apost_compare_2809_current}
    \end{subfigure}
    \caption{\label{fig:True_max_compare_current_NPBE}%
  The comparison between the maximal error estimator $\Delta_N^{\max}(\mu)$ and the true error for 
the classical NPBE during the greedy iteration at the current RB dimension $i = 1,\dots, N$ for the 
biomolecules in cases (a) to (c).}
\end{figure}

In this section, we demonstrate via the RB approximation, that the NRPBE model is more accurate and 
computationally efficient than the classical NPBE. In a similar style as in \Cref{ssec:NRPBE_RBM}, we 
consider the biomolecules in cases (a) to (c) with the corresponding domain lengths and grid dimensions. 
We demonstrate the accuracy of the RB approximation for the classical NPBE model in order to compare it with 
the NRPBE model. We begin by demonstrating in \Cref{fig:True_max_compare_current_NPBE}, the comparison of the 
error decay between the maximal error estimator $\Delta_N^{\max}(\mu)$ and the true error for the classical 
NPBE during the greedy iteration at the current RB dimension $i = 1,\dots, N$ for the biomolecules in cases 
(a) to (c). 

\setlength{\fwidth}{3.8cm}
\setlength{\fheight}{3.5cm}

We notice that the RBM constructs a small ROM (i.e., $N=2$) of high accuracy ($\mathcal{O}(10^{-12})$) 
for the NRPBE in \Cref{table:Max_vs_true_error} because of the regularized nature of the model, the RBM 
applied to the classical NPBE, nevertheless, generates a ROM of dimension $N = 6$ at the accuracy of 
$\mathcal{O}(10^{-4})$ for most biomolecules in \Cref{fig:True_max_compare_current_NPBE} 
\cite{KwBFM2017, Kweyu2016}. This is because in the latter case, the short-range component of the 
electrostatic potential impedes the reduction process due to the sharp cusps or singularities which 
are hard to capture in the ROM. Furthermore, case (a) has a slightly smaller ROM dimension due to its 
small number of atoms as compared to the rest, hence its small number of solution singularities (cusps) 
to be captured in the ROM. 

\setlength{\fwidth}{3.8cm}
\setlength{\fheight}{3.5cm}
\begin{figure}[t]
\captionsetup{width=\linewidth}
\begin{center}
    \begin{tikzpicture}
    \begin{customlegend}[legend columns=-1, legend style={/tikz/every even column/.append 
    style={column sep=2.0cm}} , legend entries={True error, $\Delta_N(\mu)$, tol}]
    \addlegendimage{blue,solid, line width = 2pt}
    \addlegendimage{black,solid, line width = 2pt}
    \addlegendimage{black,dashed, line width = 2pt}
    \end{customlegend}
    \end{tikzpicture}
    \end{center}

    \begin{subfigure}[b]{0.3\textwidth}
      \centering
      \includegraphics[height=4.2cm]{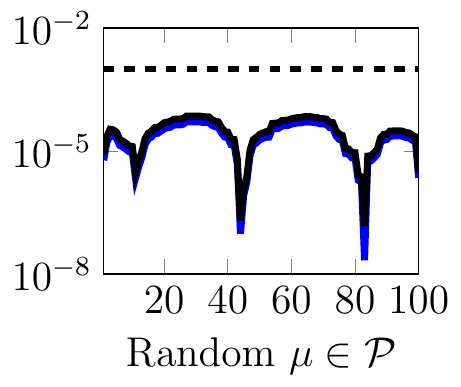}
      \caption{Error for case (a).}
      \label{fig:True_apost_compare_acet_current_rand}
    \end{subfigure}
    \hfill
    \begin{subfigure}[b]{0.3\textwidth}
        \centering
        \includegraphics[height=4.2cm]{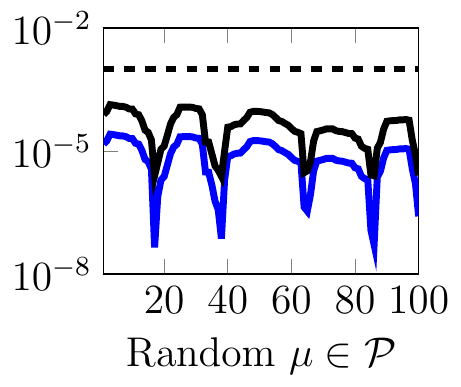}
        \caption{Error for case (b).}
        \label{fig:True_apost_compare_2809_current_rand}
    \end{subfigure}
    \hfill
    \begin{subfigure}[b]{0.3\textwidth}
        \centering
        \includegraphics[height=4.0cm]{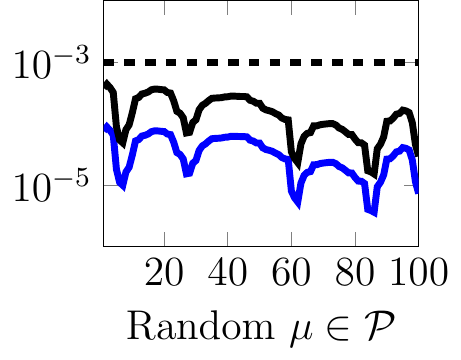}
        \caption{Error for case (c).}
        \label{fig:True_apost_compare_1FAS_final_rand}
    \end{subfigure}
    \caption{\label{fig:True_max_rand_compare_final_NPBE}%
  Comparison between the error estimator and the true error for the classical NPBE for the cases (a) to (c) 
  for the final ROM at $100$ random (varying) parameter values $\mu \in \mathcal{P}$.}
\end{figure}

The accuracy of the RB approximation of the classical NPBE in \Cref{fig:True_max_rand_compare_final_NPBE} is 
much lower than that of the NRPBE in \Cref{fig:True_max_rand_compare_final_NRPBE} due to the inaccurate 
approximation of the short-range component inherent in the former. This demonstrates the efficacy of the 
regularization scheme based on the RS tensor technique. The oscillations in the error in  
\Cref{fig:True_max_rand_compare_final_NPBE} also justifies the irregularity of the singular solution, which 
impedes the model reduction process.

\subsubsection{Runtimes and Computational Speed-ups}
\label{ssec:RPBE_RS_runtimes}

We compare the computational runtime of computing both the classical and regularized NPBE models as well 
as that of the corresponding ROM (using the RBM) in \Cref{table:Runtimes_speedups1}. The respective PBE 
models were applied to the protein \textit{fasciculin 1}. Given a fixed value of the parameter $\mu$, 
\Cref{table:Runtimes_speedups1} compares the runtimes for solving the FOM (using the FDM), constructing 
the ROM (using the RBM), and solving the ROM (using direct methods), for the classical LPBE and NPBE with 
those of the regularized LRPBE and NRPBE, respectively. It is clear that the RBM spends more time in the 
offline phase of the greedy algorithm to compute snapshots for the classical NPBE than on the NRPBE model, 
see \Cref{fig:True_max_compare_current_NPBE} and \Cref{table:Max_vs_true_error}. This is mainly because of 
the presence of rapid singularities in the PBE solution, which provides an onerous task in the construction 
of the ROM. 
\begin{table}[t]
\centering
\begin{tabular}{|c|c|c|c|c|}
  \hline
  \multicolumn{5}{|c|}{\bf{Runtime (seconds) for the PBE and the RPBE}} \\ \hline
   & LPBE & NPBE & LRPBE & NRPBE \\ \hline
  FOM & 17.68 & 34.40 & 22.83 & 28.30 \\ \hline
  RBM & 107.98 & 238.78 & 90.04 & 76.38 \\ \hline
  ROM & 2.22e-02 & 2.40e-02 & 2.10e-03 & 6.59e-03 \\ \hline
  \end{tabular}
\caption{Runtimes for the FOM, RBM and ROM for the linear and nonlinear variants of both the classical and 
the regularized PBE models.}
\label{table:Runtimes_speedups1}
\end{table}

Consequently, \Cref{table:Runtimes_speedups2}, shows that solving the FOM for the NRPBE by 
the classical numerical techniques (in this case, the FDM) is sufficient and computationally efficient 
only for a single parameter value. However, for many varying parameter values, the RBM is more 
efficient because it constructs only a small ROM once, which can then be solved fast to obtain the 
solutions at any values of the parameter. For instance, for $1000$ different parameter values, the ROM 
runtime is $\approx 6.59 \times 10^{-3} \mbox{sec}$, leading to a total runtime of $\approx 82.97 
\mbox{sec}$ to solve the NRPBE using the RBM technique instead of $\approx 28300 \mbox{sec}$ by the 
FDM solver. Note that the runtimes for the $1000$ varying parameter values for the FOM are mere 
approximations based on that of the single parameter value, since simulating the FOM for so many 
times is impractical.
  
\begin{table}[t]
\centering
\captionsetup{width=\linewidth}
\begin{tabular}{|c|c|c|c|}
\hline
  \multicolumn{4}{|c|}{\bf{Runtime (seconds) and speed-up using the FDM and the RBM}} \\
  \hline
  No. of parameters & FOM for NRPBE & RBM for NRPBE & Speed-up \\
  \hline
  1 & 28.30 & 76.38 & 0.37 \\ \hline
  10 & $\approx$ 283.00 & $\approx$ 76.44 & 3.70 \\ \hline
  100 & $\approx$ 2830.00 & $\approx$ 77.04 & 36.73 \\ \hline
  1000 & $\approx$ 28300.00 & $\approx$ 82.97 & 341.09 \\ \hline
\end{tabular}
\caption{Comparison of the runtimes and speed-ups between the FOM and the RBM for the NRPBE in a 
multiparameter context.}
\label{table:Runtimes_speedups2}
\end{table}

\section{Conclusions}\label{sec:Conclusions}

In this study, we review the salient properties of the RS canonical tensor format as a regularization 
scheme for the nonlinear PBE (NPBE) for calculating the electrostatic potential within and around 
biomolecules as proposed in \cite{BeKKKS:18,KwKKMB:18}. Among these properties is the grid-based RS tensor 
splitting of the Dirac delta distribution into the smooth and singular source term components. The NPBE is 
then discretized with the smooth approximation of the Dirac delta distribution, yielding a regularized FOM 
that is devoid of the singularities in the resultant solution. The RBM leverages this property by 
constructing a regularized ROM at extremely low computational costs as compared to that 
of the classical variant. This avoids constructing a ROM which comprises of the highly singular 
component of the electrostatic potential, thereby reducing the errors in the numerical approximation. 
The total potential is obtained by adding the regularized component (solution of the ROM), which is 
lifted (by projection) to the high-fidelity space, $\mathcal{N}$, to the directly precomputed 
canonical tensor representation of the short-range component of the Newton kernel.

\section*{Acknowledgement}
The authors thank the following organizations for financial and material support on 
this project: International Max Planck Research School (IMPRS) for Advanced Methods in Process 
and Systems Engineering and Max Planck Society for the Advancement of Science (MPG).

  \begin{footnotesize}

  \bibliographystyle{unsrtnat}
  \bibliography{BSE_Fock_Sums1.bib,Kweyu_refer.bib}

\begin{thebibliography}{52}
\providecommand{\natexlab}[1]{#1}
\providecommand{\url}[1]{\texttt{#1}}
\expandafter\ifx\csname urlstyle\endcsname\relax
  \providecommand{\doi}[1]{doi: #1}\else
  \providecommand{\doi}{doi: \begingroup \urlstyle{rm}\Url}\fi

\bibitem[Benner et~al.(2021)Benner, Khoromskaia, Khoromskij, Kweyu, and
  Stein]{BeKKKS:18}
P.~Benner, V.~Khoromskaia, B.~Khoromskij, C.~Kweyu, and M.~Stein.
\newblock Regularization of {P}oisson-{B}oltzmann type equations with singular
  source terms using the range-separated tensor format.
\newblock \emph{{SIAM} J. Sci. Comput.}, 43\penalty0 (1):\penalty0 A415--A445,
  2021.
\newblock \doi{10.1137/19M1281435}.

\bibitem[Kweyu et~al.(In prep.)Kweyu, Khoromskaia, Khoromskij, Stein, and
  Benner]{KwKKMB:18}
C.~Kweyu, V.~Khoromskaia, B.~Khoromskij, M.~Stein, and P.~Benner.
\newblock Solution decomposition for the nonlinear {P}oisson-{B}oltzmann
  equation using the range-separated tensor format.
\newblock In prep.

\bibitem[Benner et~al.(2018)Benner, Khoromskaia, and Khoromskij]{BKK_RS:18}
P.~Benner, V.~Khoromskaia, and B.~N. Khoromskij.
\newblock Range-separated tensor format for many-particle modeling.
\newblock \emph{SIAM J. Sci. Comput.}, \penalty0 (2):\penalty0 A1034--A1062,
  2018.

\bibitem[Khoromskij(2020)]{khor-DiracRS:2018}
B.~N. Khoromskij.
\newblock Range-separated tensor representation of the discretized
  multidimensional {D}irac delta and elliptic operator inverse.
\newblock \emph{J. Comp. Phys.}, 401, 2020.
\newblock \doi{10.1016/j.jcp.2019.108998}.

\bibitem[Baker(2004)]{Baker:04}
Nathan~A. Baker.
\newblock {P}oisson-{B}oltzmann methods for biomolecular electrostatics.
\newblock \emph{Methods in Enzymology}, 383:\penalty0 94--118, 2004.

\bibitem[Baker(2005)]{Baker2005}
N.~A. Baker.
\newblock \emph{Biomolecular applications of {P}oisson-{B}oltzmann equation},
  volume~21 of \emph{Reviews in Computational Chemistry}.
\newblock John Wiley \& Sons, Hoboken, NJ, USA, 2005.
\newblock \doi{10.1002/0471720895.ch5}.

\bibitem[Jason~Wagoner(2004)]{WaBa:04}
Nathan A.~Baker Jason~Wagoner.
\newblock Solvation forces on biomolecular structures: A comparison of explicit
  solvent and {P}oisson-{B}oltzmann models.
\newblock \emph{J Comput. Chem.}, 25:\penalty0 1623–1629, 2004.

\bibitem[Lu et~al.(2008)Lu, Zhou, Holst, and McCammon]{Holst:2008}
B.~Z. Lu, Y.~C. Zhou, M.~J. Holst, and J.~A. McCammon.
\newblock Recent progress in numerical methods for {P}oisson-{B}oltzmann
  equation in biophysical applications.
\newblock \emph{Commun. Comp. Phys.}, 3\penalty0 (5):\penalty0 973--1009, 2008.

\bibitem[Fogolari et~al.(2002)Fogolari, Brigo, and Molinari]{Fogolari2002}
F.~Fogolari, A.~Brigo, and H.~Molinari.
\newblock The {P}oisson-{B}oltzmann equation for biomolecular electrostatics: a
  tool for structural biology.
\newblock \emph{J. Mol. Recognit.}, 15\penalty0 (6):\penalty0 377--392, 2002.
\newblock \doi{10.1002/jmr.577}.

\bibitem[Neves-Petersen and Petersen(2003)]{Neves-Petersen2003}
M.~T. Neves-Petersen and S.~Petersen.
\newblock Protein electrostatics: {A} review of the equations and methods used
  to model electrostatic equations in biomolecules - applications in
  biotechnology.
\newblock \emph{Biotechnol. Annu. Rev.}, 9:\penalty0 315--395, 2003.
\newblock \doi{10.1016/S1387-2656(03)09010-0}.

\bibitem[Stein et~al.(2010)Stein, Gabdoulline, and Wade]{Stein:2010}
M.~Stein, R.~R. Gabdoulline, and R.~C. Wade.
\newblock Cross-species analysis of the glycoliticmpathway by comparison of
  molecular interaction fields.
\newblock \emph{Molecular Biosystems}, 6:\penalty0 162--174, 2010.

\bibitem[Holst(1994)]{Holst94}
M.~J. Holst.
\newblock \emph{Multilevel methods for the {P}oisson-{B}oltzmann equation}.
\newblock Ph.{D}. {T}hesis, Numerical Computing group, University of Illinois,
  Urbana-Champaign, IL, USA, 1994.

\bibitem[Dong et~al.(2008)Dong, Oslen, and Baker]{Dong2008}
F.~Dong, B.~Oslen, and N.~A. Baker.
\newblock Computational methods for biomolecular electrostatics.
\newblock \emph{Methods Cell Biol}, 84\penalty0 (1):\penalty0 843--870, 2008.
\newblock \doi{10.1016/S0091-679X(07)84026-X}.

\bibitem[Baker et~al.(2001{\natexlab{a}})Baker, Holst, and Wang]{Baker2001}
N.~A. Baker, M.~J. Holst, and F.~Wang.
\newblock The adaptive multilevel finite element solution of the
  {P}oisson-{B}oltzmann equation on massively parallel computers.
\newblock \emph{IBM J. Res. Devel.}, 45:\penalty0 427--438, 2001{\natexlab{a}}.

\bibitem[Wang and Luo(2010)]{Wang2010}
J.~Wang and R.~Luo.
\newblock Assessment of linear finite difference {P}oisson-{B}oltzmann solvers.
\newblock \emph{J. Comput. Chem.}, 31:\penalty0 1689--1698, 2010.
\newblock \doi{10.1016/j.cpc.2015.08.029}.

\bibitem[Holst et~al.(2000)Holst, Baker, and Wang]{Holst2000}
M.~Holst, N.~Baker, and F.~Wang.
\newblock Adaptive multilevel finite element solution of the
  {P}oisson-{B}oltzmann equation: algorithms and examples.
\newblock \emph{J. Comp. Chem.}, 21:\penalty0 1319--1342, 2000.
\newblock \doi{10.1002/1096-987X(20001130)21:15<1319::AID-JCC1>3.0.CO;2-8}.

\bibitem[Boschitsch and Fenley(2004)]{Boschitsch2004}
A.~H. Boschitsch and M.~O. Fenley.
\newblock Hybrid boundary element and finite difference method for solving the
  nonlinear {P}oisson-{B}oltzman equation.
\newblock \emph{J. Comput. Chem.}, 25\penalty0 (7):\penalty0 935--955, 2004.
\newblock \doi{10.1002/jcc.20000}.

\bibitem[Zhou(1993)]{Zhou1993}
H.~X. Zhou.
\newblock Boundary element solution of macromolecular electrostatics:
  {I}nteaction energy between two proteins.
\newblock \emph{Biophys. J.}, 65\penalty0 (2):\penalty0 955--963, 1993.
\newblock \doi{10.1016/S0006-3495(93)81094-4}.

\bibitem[Xie(2014)]{Xie:14}
D.~Xie.
\newblock New solution decomposition and minimization scheme for
  {P}oisson-{B}oltzmann equation in calculation of biomolecular electrostatics.
\newblock \emph{J Comp. Phys.}, 275:\penalty0 294--309, 2014.

\bibitem[Chen et~al.(2007)Chen, Holst, and Xu]{Chen:07}
L.~Chen, M.J. Holst, and J.~Xu.
\newblock The finite element approximation of the nonlinear
  {P}oisson-{B}oltzmann equation.
\newblock \emph{{SIAM} J. Numer. Anal.}, 45\penalty0 (6):\penalty0 2298--2320,
  2007.
\newblock \doi{10.1137/060675514}.

\bibitem[Mirzadeh et~al.(2013)Mirzadeh, Theillard, Helgadottir, Boy, and
  Gibou]{Mirzadeh:13}
M.~Mirzadeh, M.~Theillard, A.~Helgadottir, D.~Boy, and F.~Gibou.
\newblock An adaptive, finite difference solver for the nonlinear
  {P}oisson-{B}oltzmann equation with applications to biomolecular
  computations.
\newblock \emph{Commun. Comput. Phys.}, 13\penalty0 (1):\penalty0 150--173,
  2013.
\newblock \doi{10.4208/cicp.290711.181011s}.

\bibitem[Ji et~al.(2019)Ji, Chen, and Xu]{Ji:2018}
L.~Ji, Y.~Chen, and Z.~Xu.
\newblock A reduced basis method for the nonlinear {P}oisson-{B}oltzmann
  equation.
\newblock \emph{Adv. Appl. Math. Mech.}, 11:\penalty0 1200--1218, 2019.
\newblock \doi{10.4208/aamm.OA-2018-0188}.

\bibitem[Chaturantabut and Sorensen(2010)]{morChaS10}
S.~Chaturantabut and D.~C. Sorensen.
\newblock Nonlinear model reduction via discrete empirical interpolation.
\newblock \emph{{SIAM} J. Sci. Comput.}, 32\penalty0 (5):\penalty0 2737--2764,
  2010.
\newblock \doi{10.1137/090766498}.

\bibitem[Grepl et~al.(2007)Grepl, Maday, Nguyen, and Patera]{morGreMNetal07}
M.~A. Grepl, Y.~Maday, N.~C. Nguyen, and A.~T. Patera.
\newblock Efficient reduced-basis treatment of nonaffine and nonlinear partial
  differential equations.
\newblock \emph{{ESAIM}: Math. Model. Numer. Anal.}, 41\penalty0 (3):\penalty0
  575--605, 2007.

\bibitem[Barrault et~al.(2004)Barrault, Maday, Nguyen, and
  Patera]{morBarMNetal04}
M.~Barrault, Y.~Maday, N.~C. Nguyen, and A.~T. Patera.
\newblock An `empirical interpolation' method: application to efficient
  reduced-basis discretization of partial differential equations.
\newblock \emph{C. R. Math. Acad. Sci. Paris}, 339\penalty0 (9):\penalty0
  667--672, 2004.
\newblock ISSN 1631-073X.

\bibitem[Sharp and Honig(1990)]{SharpHonig90}
K.~A. Sharp and B.~Honig.
\newblock Electrostatic interactions in macromolecules: theory and
  applications.
\newblock \emph{Annu. Rev. Biophys. Chem.}, 19:\penalty0 301--332, 1990.

\bibitem[Kweyu et~al.(2020)Kweyu, Feng, Stein, and Benner]{KwBFM2017}
C.~Kweyu, L.~Feng, M.~Stein, and P.~Benner.
\newblock Fast solution of the {P}oisson-{B}oltzmann equation with nonaffine
  parametrized boundary conditions using the reduced basis method.
\newblock \emph{Computing and Visualization in Science}, 23\penalty0 (15),
  2020.
\newblock \doi{10.1007/s00791-020-00336-z}.

\bibitem[Fogolari et~al.(1999)Fogolari, Zuccato, Esposito, and
  Viglino]{Fogolari99}
F.~Fogolari, P.~Zuccato, G.~Esposito, and P.~Viglino.
\newblock Biomolecular electrostatics with the linearized {P}oisson-{B}oltzmann
  equation.
\newblock \emph{Biophys. J.}, 76\penalty0 (1):\penalty0 1--16, 1999.
\newblock \doi{10.1016/S0006-3495(99)77173-0}.

\bibitem[Qin et~al.(2010)Qin, Meng-Juei, Jun, and Ray]{Qin_Ray_2010}
C.~Qin, H.~Meng-Juei, W.~Jun, and L.~Ray.
\newblock Performance of nonlinear finite-difference {P}oisson-{B}oltzmann
  solvers.
\newblock \emph{Journal of Chemical Theory and Computation}, 6\penalty0
  (1):\penalty0 203--211, 2010.
\newblock \doi{10.1021/ct900381r}.
\newblock URL \url{https://doi.org/10.1021/ct900381r}.

\bibitem[Chern et~al.(2003)Chern, Liu, and Wang]{Chern:2003}
I.~Chern, J.~Liu, and W.~Wang.
\newblock Accurate evaluation of electrostatics for macromolecules in solution.
\newblock \emph{Methods Appl. Anal.}, 10\penalty0 (2):\penalty0 309--328, 2003.

\bibitem[of~Mathematics()]{Newt_kernel}
Encyclopedia of~Mathematics.
\newblock Newton potential.
\newblock
  \url{http://www.encyclopediaofmath.org/index.php?title=Newton_potential&oldid=33114}.
\newblock Accessed: 2018-03-12.

\bibitem[Khoromskaia and Khoromskij(2014)]{VeBoKh:Ewald:14}
V.~Khoromskaia and B.~N. Khoromskij.
\newblock Grid-based lattice summation of electrostatic potentials by assembled
  rank-structured tensor approximation.
\newblock \emph{Comp. Phys. Comm.}, 185\penalty0 (12), 2014.

\bibitem[Khoromskij and Khoromskaia(2009)]{khor-ml-2009}
B.~N. Khoromskij and V.~Khoromskaia.
\newblock Multigrid accelerated tensor approximation of function related
  multidimensional arrays.
\newblock \emph{SIAM J. Sci. Comp.}, 31\penalty0 (4):\penalty0 3002--3026,
  2009.
\newblock \doi{10.1137/080730408}.

\bibitem[Rocchia et~al.(2001)Rocchia, Alexov, and Honig]{Rocchia_2001}
W.~Rocchia, E.~Alexov, and B.~Honig.
\newblock Extending the applicability of the nonlinear {P}oisson-{B}oltzmann
  equation: multiple dielectric constants and multivalent ions.
\newblock \emph{J. Phys. Chem.}, 105\penalty0 (28):\penalty0 6507--6514, 2001.
\newblock \doi{10.1021/jp010454y}.

\bibitem[Luty et~al.(1992)Luty, Davis, and McCammon]{Brock_1992}
B.A. Luty, M.E. Davis, and J.A McCammon.
\newblock Solving the finite-difference nonlinear {P}oisson-{B}oltzmann
  equation.
\newblock \emph{J. Comput. Chem.}, 13\penalty0 (9):\penalty0 1114--1118, 1992.
\newblock \doi{10.1002/jcc.540130911}.

\bibitem[Oberoi and Allewell(1993)]{Oberoi_1993}
H.~Oberoi and N.~M. Allewell.
\newblock Multigrid solution of the nonlinear {P}oisson-{B}oltzmann equation
  and calculation of titration curves.
\newblock \emph{Biophys. J.}, 65\penalty0 (1):\penalty0 48--55, 1993.
\newblock \doi{10.1016/S0006-3495(93)81032-4}.

\bibitem[Holst and Saied(1995)]{Holst:95}
M.~Holst and F.~Saied.
\newblock Numerical solution of the nonlinear {P}oisson-{B}oltzmann equation:
  Developing more robust and efficient methods.
\newblock \emph{J. Comput. Chem.}, 16:\penalty0 337--364, 1995.

\bibitem[Shestakov et~al.(2002)Shestakov, Milovich, and Noy]{Shestakov:2002}
A.~I. Shestakov, J.~L. Milovich, and A.~Noy.
\newblock Solution of the nonlinear {P}oisson-{B}oltzmann equation using
  pseudo-transient continuation and the finite element method.
\newblock \emph{Commun. Comput. Phys.}, 247:\penalty0 62--79, 2002.
\newblock \doi{10.1006/jcis.2001.8033}.

\bibitem[Notay(2010)]{Notay:2010}
Y.~Notay.
\newblock An aggregation-based algebraic multigrid method.
\newblock \emph{Electronic Transactions on Numerical Analysis}, 37:\penalty0
  123--146, 2010.

\bibitem[Hesthaven et~al.(2016)Hesthaven, Rozza, and Stamm]{morHest16}
J.S. Hesthaven, G.~Rozza, and B.~Stamm.
\newblock \emph{Certified {R}educed {B}asis {M}ethods for {P}arametrized
  {P}artial {D}ifferential {E}quations}.
\newblock Springer International Publishing, 2016.
\newblock \doi{10.1007/978-3-319-22470-1}.

\bibitem[Benner et~al.(2015)Benner, Gugercin, and Willcox]{morBenGW15}
P.~Benner, S.~Gugercin, and K.~Willcox.
\newblock A survey of model reduction methods for parametric systems.
\newblock \emph{SIAM Review}, 57\penalty0 (4):\penalty0 483--531, 2015.
\newblock \doi{10.1137/130932715}.

\bibitem[Eftang(2011)]{morEft11}
J.~L. Eftang.
\newblock \emph{Reduced basis methods for parametrized partial differential
  equations}.
\newblock Ph.{D}. {T}hesis, Norwegian University of Science and Technology,
  Trondheim, Norway, 2011.

\bibitem[Rozza et~al.(2008)Rozza, Huynh, and Patera]{morRozHP08}
G.~Rozza, D.~B.~P. Huynh, and A.~T. Patera.
\newblock Reduced basis approximation and a posteriori error estimation for
  affinely parametrized elliptic coercive partial differential equations.
\newblock \emph{Archives of Computational Methods in Engineering}, 15\penalty0
  (3):\penalty0 229--275, 2008.
\newblock \doi{10.1007/s11831-008-9019-9}.

\bibitem[Volkwein(2013)]{morVol13}
S.~Volkwein.
\newblock Model reduction using proper orthogonal decomposition.
\newblock Lecture notes, University of Konstanz, 2013.

\bibitem[Feng and Benner(2014)]{morBenF14}
L.~Feng and P.~Benner.
\newblock \emph{Reduced Order Methods for modeling and computational reduction,
  MS\&A Series}, volume~9, chapter 6: A robust algorithm for parametric model
  order reduction based on implicit moment matching, pages 159--186.
\newblock Springer-Verlag, Berlin, Heidelberg, New York, 2014.
\newblock \doi{10.1007/978-3-319-02090-7_6}.

\bibitem[Feng et~al.(2017)Feng, Mangold, and Benner]{Feng2016}
L.~Feng, M.~Mangold, and P.~Benner.
\newblock Adaptive {POD-DEIM} basis construction and its application to a
  nonlinear population balance system.
\newblock \emph{AIChE Journal}, pages 3832--3844, 2017.
\newblock \doi{10.1002/aic.15749}.

\bibitem[Baker et~al.(2001{\natexlab{b}})Baker, Sept, Joseph, Holst, and
  McCammon]{Bakersept2001}
N.~A. Baker, D.~Sept, S.~Joseph, M.~J. Holst, and J.~A. McCammon.
\newblock Electrostatics of nanosystems: application to microtubules and the
  ribosome.
\newblock \emph{Proc. Nat. Acad. Sci. U.S.A.}, 98\penalty0 (18):\penalty0
  10037--10041, 2001{\natexlab{b}}.
\newblock \doi{10.1073/pnas.181342398}.

\bibitem[Holst and Saied(1993)]{Holst:93}
M.~Holst and F.~Saied.
\newblock Multigrid solution of the {P}oisson-{B}oltzmann equation.
\newblock \emph{J. Comput. Chem.}, 14:\penalty0 105--113, 1993.

\bibitem[Vergara-Perez and Marucho(2016)]{Vergara-Perez2016}
S.~Vergara-Perez and M.~Marucho.
\newblock {MPBEC}, a {M}atlab program for biomolecular electrostatic
  calculations.
\newblock \emph{Comput. Phys. Commun.}, 198:\penalty0 179--194, 2016.
\newblock \doi{10.1016/j.cpc.2015.08.029}.

\bibitem[le~Du et~al.(1992)le~Du, Marchot, Bougis, and
  Fontecilla-Camps]{DuMaBoFo:92}
M.H. le~Du, P.~Marchot, P.E. Bougis, and J.C. Fontecilla-Camps.
\newblock 1.9 {A}ngstrom resolution structure of fasciculine 1, an
  anti-acetylcholinesterase toxin from green mamba snake venom.
\newblock \emph{J. Biol. Chem.}, 267:\penalty0 22122--22130, 1992.

\bibitem[Hinds et~al.(1998)Hinds, Maurer, J., and Nicola]{HinMauZhaNic:98}
M.G. Hinds, T.~Maurer, Zhang J., and N.A. Nicola.
\newblock Solution structure of {L}eukemia inhibitory factor.
\newblock \emph{BiolChem}, 273:\penalty0 13738--13745, 1998.
\newblock \doi{10.1074/jbc.273.22.13738}.

\bibitem[Kweyu et~al.(2016)Kweyu, Hess, Feng, Stein, and Benner]{Kweyu2016}
C.~Kweyu, M.~Hess, L.~Feng, M.~Stein, and P.~Benner.
\newblock Reduced basis method for {P}oisson-{B}oltzmann {E}quation.
\newblock In M.~Papadrakakis, V.~Papadopoulos, G.~Stefanou, and V.~Plevris,
  editors, \emph{ECCOMAS Congress 2016 - Proc. of the {VII} {E}uropean Congress
  on Computational Methods in Applied Sciences and Engineering}, volume~2,
  pages 4187--4195, Athens, 2016. National Technical University of Athens.

\end{thebibliography}

  \end{footnotesize}

\end{document}